\documentclass[OptSoft]{jtcam_preprint}

\usepackage[T1]{fontenc}
\usepackage{float}
\usepackage{algorithm}
\usepackage{amsmath}
\usepackage{graphicx}
\providecommand{\dif}{\mathrm{d}}
\usepackage{graphicx,wasysym}
\usepackage[T1]{fontenc}
\usepackage{calrsfs}
\usepackage{nicefrac}
\usepackage{babel}

\DeclareRobustCommand{\PhaseOne}{[\tikz[baseline=-\the\dimexpr\fontdimen22\textfont2\relax,inner sep=0pt,thick]\draw[red](0,0)--(.5,0);]}
\DeclareRobustCommand{\PhaseTwo}{[\tikz[baseline=-\the\dimexpr\fontdimen22\textfont2\relax,inner sep=0pt,thick]\draw[blue](0,0)--(.5,0);]}

\title{An implicit staggered algorithm for CPFEM-based analysis
of aluminum}
\runningtitle{Implicit staggered algorithm for CPFEM}

\author[1,2]{P. Areias}
\author[3]{Charles dos Santos}
\author[2,4]{R. Melicio}
\author[1,2]{N. Silvestre}
\runningauthor{P. Areias et al.}
\affil[1]{IST - Departament of Mechanical Engineering, Instituto Superior T\'ecnico\\ Universidade de Lisboa\\ Av. Rovisco Pais 1, 1049-001, Lisboa, Portugal} 
\affil[2]{IDMEC - Instituto de Engenharia Mec\^{a}nica\\ Av. Rovisco Pais 1, 1049-001, Lisboa, Portugal} 
\affil[3]{INSA - Institut National des Sciences Appliqu\'ees (alumni) \\ Rouen Normandie, France} 
\affil[4]{AEROG - Universidade da Beira Interior\\ Cal\c{c}ada Fonte do Lameiro, 6201-001, Covilh\~{a}, Portugal} 


\begin{document}


\nocite{*} 	
\begin{abstract}
In this paper, we propose an implicit staggered algorithm for crystal
plasticity finite element method (CPFEM) which makes use of dynamic
relaxation at the constitutive integration level. An uncoupled version
of the constitutive system consists of a multi-surface flow law complemented
by an evolution law for the hardening variables. Since a saturation
law is adopted for hardening, a sequence of nonlinear iteration followed
by a linear system is feasible. To tie the constitutive unknowns,
the dynamic relaxation method is adopted. A Green-Naghdi plasticity
model is adopted based on the Hencky strain calculated using a $[2/2]$
Pad\'e approximation. For the incompressible case, the approximation
error is calculated exactly. A enhanced-assumed strain (EAS) element
technology is adopted, which was found to be especially suited to
localization problems such as the ones resulting from crystal plasticity
plane slipping. Analysis of the results shows significant reduction
of drift and well defined localization without spurious modes or hourglassing. 
\end{abstract}
\section{Introduction}

Realistic simulations of anisotropy in metal polycrystals require
robust single crystal algorithms that consistently produce results
within an established error bound. In addition, since polycrystals
typically include a large number of grains and significant computational
costs, efficiency improvements are a necessity. Fully coupled constitutive
systems for single crystals involve a large number of constitutive
unknowns (plastic strain tensor, hardening variables, among others)
in a nonlinear and frequently nonsmooth system. In a single FCC simulation
at one quadrature point, there are at least $18$ constitutive unknowns,
corresponding to hardening in each of the $12$ dominant slip systems
plus $6$ unknowns corresponding to the flow law (either plastic strain
or final stress). Computational costs of the fully coupled constitutive
system in a polycrystal are prohibitive for practical applications
without considerable computational investment. Staggered algorithms
are often used, but naive implementations produce drifting, since
there is no measure to reduce the solution error. In that case, both
time step and sequence of operations affect the solution \cite{matthies2006}.
A comprehensive description is given by Felippa, Park, and Farhat
\cite{felippa2001}. In terms of solution classes for the coupled
equations, options are as follows:
\begin{itemize}
\item Classical full Newton iteration for the coupled constitutive system
(e.g. \cite{neto2008}).
\item Sequential solution of the constitutive equations, in a staggered
form.
\item Fixed point iteration with relaxation (e.g. in the fluid-structure
interaction context \cite{letallec2001,kuttler2008}). 
\item Block-Newton partition approach with approximate out-of-diagonal blocks
(e.g. \cite{matthies2003}).
\item Block-Newton-Krylov where only multiplications are used for the out-of-diagonal
blocks \cite{michler2006}.
\end{itemize}
Effective implicit staggered algorithms have been developed for thermoelasticity
\cite{erbts2012} and phase-field simulations \cite{schapira2023}.
These are here adopted to improve the accuracy of the staggered algorithm.
For FCC (aluminum), we propose an implicit staggered algorithm which
strongly couples the solutions (in the sense of Matthies, Niekamp
and Steindorf \cite{matthies2006}) to avoid drifting. 

Classical works on single crystal plasticity emphasize the strain
localization physics and describe the essentials of what is now CPFEM,
see \cite{rice1971,asaro1977}. Significant developments were achieved
to incorporate mettallurgical effects in hardening and coupling with
grain boundary phenomena (see, e.g. \cite{cailletaud2003}). Current
success of algorithms for single crystals can be observed by large-scale
polycrystalline ensembles, either combined with homogenization or
not. Currently, theoretically sound frameworks exist for single crystal
plasticity, including that of Kaiser and Menzel \cite{kaiser2019},
where gradient effects in hardening are considered. A review of developments
in single crystal plasticity up to the year 2010 was performed by
Roters \emph{et al.} \cite{roters2010}. 

This work is organized as follows: in section \Cref{sec:Aluminum-crystal-plasticity},
a description of the aluminum FCC crystal plasticity is presented,
in Section \Cref{sec:Constitutive-integration} the constitutive integration
algorithm for the flow law and the hardening evolution law is presented
and a localization test is performed. Section \Cref{secimplstagg}
presents the proposed implicit staggered algorithm based on dynamic
relaxation, as well as a verification test for its effectiveness.
Section \Cref{secfiniteelementformulation} presents the finite
element technology, specifically a 3D enhanced assumed strain hexahedron
(EAS) which is able to capture strain localization. A polycrystal
numerical test, following the data of Alankar, Mastorakos and Field
\cite{alankar2009}, is presented in section \Cref{secnumass}.
In section \Cref{seconclusion}, conclusions are drawn with respect
to the proposed algorithm.
\section{Aluminum crystal plasticity}

\label{sec:Aluminum-crystal-plasticity}

Significant literature exists concerning aluminum plasticity, both
in phenomenological \cite{barlat1991,barlat2005} as well as single
crystal \cite{alankar2009,kasemer2020,romanova2022} cases. The FCC
dominant slip systems (which is considered here for pure Aluminum)
consist of $\{111\}$ planes and $<110>$ directions. Table \ref{tab:planes}
presents the dominant slip systems. The lattice elasticity matrix
is anisotropic and given in sub-table \ref{tab:Lattice-elasticity-matrix}.
\begin{center}
\begin{table}[h]
\begin{centering}
\caption{\label{tab:planes}Dominant slip systems of an FCC single crystal}
\par\end{centering}
\centering{}%
\begin{tabular}{ccc}
Slip system & Dense plane ($\boldsymbol{m}_{\alpha}$) & Dense direction ($\boldsymbol{n}_{\alpha}$)\tabularnewline

$1$ & $(111)$ & $[01\overline{1}]$\tabularnewline
$2$ & $(111)$ & $[10\overline{1}]$\tabularnewline
$3$ & $(111)$ & $[1\overline{1}0]$\tabularnewline
$4$ & $(11\overline{1})$ & $[011]$\tabularnewline
$5$ & $(11\overline{1})$ & $[\overline{1}0\overline{1}]$\tabularnewline
$6$ & $(11\overline{1})$ & $[\overline{1}10]$\tabularnewline
$7$ & $(1\overline{1}1)$ & $[011]$\tabularnewline
$8$ & $(1\overline{1}1)$ & $[\overline{1}01]$\tabularnewline
$9$ & $(1\overline{1}1)$ & $[\overline{1}\overline{1}0]$\tabularnewline
$10$ & $(\bar{1}11)$ & $[01\overline{1}]$\tabularnewline
$11$ & $(\bar{1}11)$ & $[\overline{1}0\overline{1}]$\tabularnewline
$12$ & $(\bar{1}11)$ & $[\overline{1}\overline{1}0]$\tabularnewline
\end{tabular}
\end{table}
\par\end{center}

\begin{table}
\begin{centering}
\subfloat[\label{tab:Lattice-elasticity-matrix}Lattice elasticity matrix]{
\begin{centering}
\vspace{0.2cm}
\par\end{centering}
\centering{}$\text{\ensuremath{\mathcal{C}=}}\left[\begin{array}{c|c|c|c|c|c}
106.75 & 60.41 & 60.41 & 0 & 0 & 0\\
\hline 60.41 & 106.75 & 60.41 & 0 & 0 & 0\\
\hline 60.41 & 60.41 & 106.75 & 0 & 0 & 0\\
\hline 0 & 0 & 0 & 28.34 & 0 & 0\\
\hline 0 & 0 & 0 & 0 & 28.34 & 0\\
\hline 0 & 0 & 0 & 0 & 0 & 28.34
\end{array}\right]\times10^{9}\quad\textrm{Pa}$}
\par\end{centering}
\begin{centering}
\subfloat[\label{tab:Hardening-properties}Hardening properties]{

\vspace{0.2cm}
\centering{}%
\begin{tabular}{c|c|c|c|c|c|c}
$\dot{\gamma}_{0}$ $\left[\textrm{s}^{-1}\right]$ & $h_{0}$ $\left[\textrm{Pa}\right]$ & $\xi_{0}$ $\left[\textrm{Pa}\right]$ & $\xi_{\infty}$ $\left[\textrm{Pa}\right]$ & $\xi_{\infty}^{\star}$ $\left[\textrm{Pa}\right]$ & $q$ $\left[-\right]$ & $n$\tabularnewline
\hline 
$0.001$ & $75\times10^{6}$ & $31\times10^{6}$ & $63\times10^{6}$ & $7\times10^{6}$ & $1.4$ & $30$\tabularnewline
\end{tabular}}
\par\end{centering}
\caption{Relevant properties for single crystal plasticity of aluminum (see
\cite{alankar2009}). These are assumed constant for $\alpha=1,\ldots,12$. }
\end{table}

\begin{figure}
\begin{centering}
\includegraphics[width=8cm]{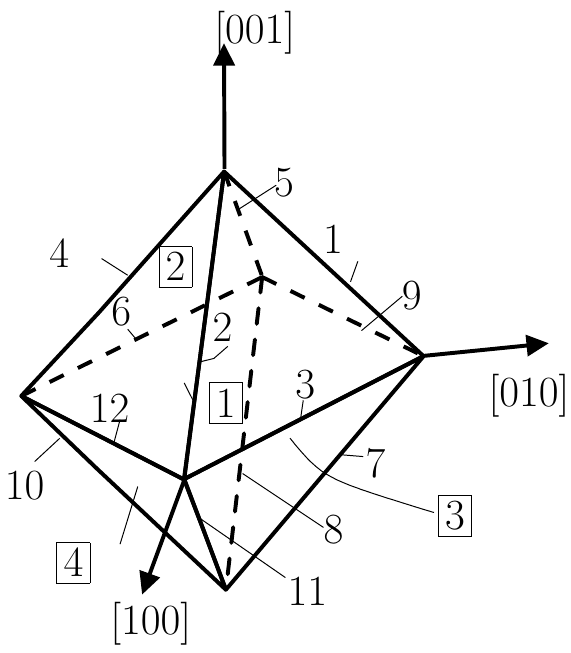}
\par\end{centering}
\caption{\label{figaDominant-slip-systems}Dominant slip systems for a FCC
crystal.}

\end{figure}

\begin{figure}
\begin{centering}
\includegraphics[width=12cm]{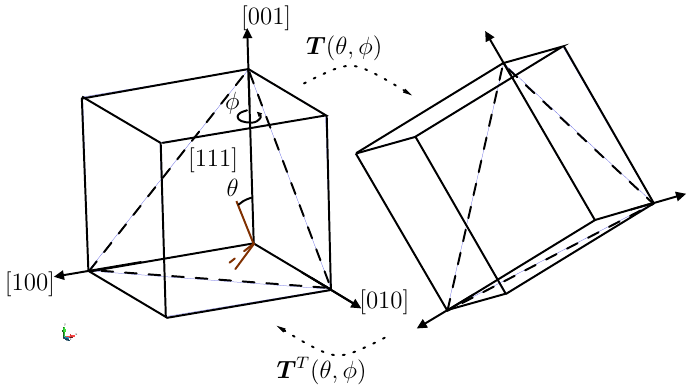}
\par\end{centering}
\caption{\label{figaRotation-of-the}Rotation of the FCC cell in the space
$x,y,z$. }

\end{figure}

Planes and directions are rotated for the analysis, as Figure \ref{figaRotation-of-the}
shows. The crystal orientation is defined in spherical coordinates
using two angles: $\theta$ and $\phi$. Figure \ref{figaRotation-of-the}
exhibits the transformation. The corresponding transformation matrix
defined as follows:

\begin{equation}
\boldsymbol{T}\left(\theta,\phi\right)=\left[\begin{array}{c|c|c}
\cos\theta\cos\phi & \cos\theta\sin\phi & -\sin\theta\\
\hline -\sin\phi & \cos\phi & 0\\
\hline \sin\theta\cos\phi & \sin\theta\sin\phi & \cos\theta
\end{array}\right]\label{eq:tthetaphi}
\end{equation}

For the stereographic representation of the loading direction, the
transpose of $\boldsymbol{T\left(\theta,\phi\right)}$ is adopted.

We use a visco-plastic formulation, which is comprised of the following
ingredients:
\begin{itemize}
\item Lattice elastic law.
\item Additive decomposition of the logarithmic strain into elastic and
plastic terms.
\item Schmidt flow law.
\item Hardening laws for the $12$ systems.
\end{itemize}
We adopt the Lagrangian Hencky strain, see \cite{schroder2002}:

\begin{equation}
\boldsymbol{\varepsilon}=\frac{1}{2}\log\left[2\boldsymbol{E}+\boldsymbol{I}\right]\label{eq:log}
\end{equation}
\\
where $\boldsymbol{E}$ is the Green--Lagrange strain and $\boldsymbol{I}$ is the
identity matrix. Schr\"{o}der, Gruttmann, and L\"{o}blein \cite{schroder2002}
as well as Shutov and co-workers \cite{shutov2014} and Miehe's group
\cite{miehe2002} have explored the additive decomposition of $\boldsymbol{\varepsilon}$
and the corresponding conjugate stress. For metals, it is a long-standing
procedure to adopt a Hooke-like law with the strain (\ref{eq:log}),
see \cite{anand1979}, which is compatible with hyperelasticity if
the Kirchhoff stress is used. The present approach is aligned with
the spirit of Asaro in \cite{asaro1983} for the constitutive law.
A review is presented by Xiao \cite{xiao2005}. Computational origins
of the additive strain decomposition in finite strains were established
by Papadopoulus and Lu \cite{papadopoulos1998,papadopoulos2001}.
The approach inherits some of the formalism adopted in small strain
elasto-plasticity. In particular, the additive decomposition, into
elastic ($\boldsymbol{\varepsilon}_{e}$) and plastic ($\boldsymbol{\varepsilon}_{p}$)
parts is retained:

\begin{equation}
\boldsymbol{\varepsilon}=\boldsymbol{\varepsilon}_{e}+\boldsymbol{\varepsilon}_{p}\label{eq:decompstess}
\end{equation}

Theoretical foundations for the use of (\ref{eq:decompstess}) in
finite strains were established by Green and Naghdi \cite{green1965}
and have been followed by schools of continuum mechanics, see Lehmann
\cite{lehmann1991}. The lattice elastic law is written as follows:

\begin{equation}
\boldsymbol{\varepsilon}_{e}=\mathcal{C}^{-1}:\boldsymbol{\sigma}\label{eq:sigmaeps}
\end{equation}

In the isotropic case, Bruhns, Xiao and Meyers \cite{bruhns2001}
proved the satisfaction of the Legendre-Hadamard condition if every
principal stretch satisfies $\lambda_{k}\in[0.21162...,1.39561...]$
which encompasses the range of elastic behavior of metals discussed
by Anand \cite{anand1979}. In (\ref{eq:decompstess}), the total
strain $\boldsymbol{\varepsilon}$ is assumed to be known from $\boldsymbol{E}$.
That allows to establish, in equation (\ref{eq:decompstess}), a coupling
between the macroscopic stress $\boldsymbol{\sigma}$ and the plastic strain
$\boldsymbol{\varepsilon}_{p}$. The stress for the system $\alpha$ is obtained
using Cauchy's lemma:

\begin{align}
\tau^{\alpha} & =\boldsymbol{m}_{\alpha}\cdot\left(\boldsymbol{\sigma}\cdot\boldsymbol{n}_{\alpha}\right)\Leftrightarrow\label{eq:tauresolved}\\
\tau^{\alpha} & =\boldsymbol{M}_{\alpha}:\boldsymbol{\sigma}\\
\tau^{\alpha} & =\boldsymbol{P}_{\alpha}:\left(\boldsymbol{\boldsymbol{\varepsilon}}-\boldsymbol{\boldsymbol{\varepsilon}}_{p}\right)
\end{align}
\\
where $\boldsymbol{P}_{\alpha}=\left(\boldsymbol{m}_{\alpha}\otimes\boldsymbol{n}_{\alpha}\right):\mathcal{C}$.
The strain rate corresponding to each dominant slip system $\alpha=1,\ldots,12$
is introduced as $\dot{\gamma}_{\alpha}$. The single crystal flow
law is provided by the Schmid relation, which is known to be an acceptable
starting point for FCC crystals:

\begin{equation}
\dot{\boldsymbol{\varepsilon}}_{p}=\sum_{\alpha=1}^{12}\boldsymbol{M}_{\alpha}\dot{\gamma}^{\alpha}\label{eq:flowlaw}
\end{equation}
\\
where the flow vector $\boldsymbol{M}_{\alpha}$ is obtained as $\boldsymbol{M}_{\alpha}=\left[\boldsymbol{m}_{\alpha}\otimes\boldsymbol{n}_{\alpha}\right]_{\textrm{symm}}$.
We note that $\boldsymbol{M}_{\alpha}$ is established for a given crystalline
structure, which is fixed during the analysis. Since symmetry ensures
that a Voigt form can be established for $\boldsymbol{M}_{\alpha}$, it becomes
possible to write: 

\begin{equation}
\boldsymbol{P}_{\alpha}=\boldsymbol{M}_{\alpha}:\mathcal{C}\label{eq:mac}
\end{equation}

If a fixed yield stress is assumed, then

\begin{equation}
\xi^{\alpha}-|\tau^{\alpha}|\geq0\label{eq:xitau}
\end{equation}

Note that, in (\ref{eq:flowlaw}), the equivalent strain rate in system
$\alpha$, $\dot{\gamma}^{\alpha}$ can assume either negative or
positive values. It is often assumed that $\dot{\gamma}^{\alpha}$
follows a visco-plastic law:

\begin{equation}
\dot{\gamma}^{\alpha}=\dot{\gamma}_{0}^{\alpha}\left|\frac{\tau^{\alpha}}{\xi^{\alpha}}\right|^{n}\textrm{sign}\left[\tau^{\alpha}\right]\label{eq:gammadotalpha}
\end{equation}

The hardening law for the critically resolved shear stress $\xi^{\alpha}$
follows Kasemer \cite{kasemer2020} saturation proposal:

\begin{equation}
\dot{\xi}^{\alpha}=h_{0}\sum_{\beta=1}^{12}\left[\left|\dot{\gamma}^{\beta}\right|\left(1-\frac{\xi^{\beta}}{\xi_{\infty}^{\beta}}\right)h_{\alpha\beta}\right]\label{eq:hard}
\end{equation}
where the coupling matrix is given by \cite{chang1981} with $q=1.4$
\cite{kasemer2020}, see also \cite{bassani1991,zhang2016}.

\begin{equation}
h_{\alpha\beta}=\begin{cases}
1 & \alpha=\beta\\
1.4 & \alpha\neq\beta
\end{cases}\label{eq:hab}
\end{equation}

Initial conditions for the previous constitutive system are:

\begin{align}
\left.\boldsymbol{\varepsilon}_{p}\right|_{0} & =\boldsymbol{0}\label{eq:epzero}\\
\xi^{\alpha} & =\xi_{0}\label{eq:xi0}
\end{align}

\section{Constitutive integration and testing}

\label{sec:Constitutive-integration}

Using the backward-Euler integration method between time steps $s$
and $s+1$, the flow law (\ref{eq:flowlaw}) reads:

\begin{equation}
\boldsymbol{\varepsilon}_{s+1}^{p}=\boldsymbol{\varepsilon}_{s}^{p}+\sum_{\alpha=1}^{12}\boldsymbol{M}_{\alpha}\Delta\gamma_{s+1}^{\alpha}\label{eq:ep}
\end{equation}

Given (\ref{eq:ep}), the \emph{trial} shear stress for system $\alpha$
is given by:

\begin{equation}
\tau_{\alpha}^{\star}=\boldsymbol{P}_{\alpha}\cdot\left(\boldsymbol{\varepsilon}-\boldsymbol{\varepsilon}_{s}^{p}\right)\label{eq:taua}
\end{equation}
\\
Trial $\tau_{\alpha}^{\star}$ is such that $\tau_{\alpha}^{\star}\Delta\gamma_{s+1}^{\alpha}\geq0$,
without sum on $\alpha$. Integration of (\ref{eq:ep}) for a given
time increment $\Delta t$ is as follows:

\begin{equation}
\boldsymbol{\varepsilon}_{s+1}^{p}-\boldsymbol{\varepsilon}_{s}^{p}-\Delta t\sum_{\alpha=1}^{12}\boldsymbol{M}_{\alpha}\dot{\gamma}_{0}^{\alpha}\left|\frac{\boldsymbol{P}_{\alpha}:\left(\boldsymbol{\boldsymbol{\varepsilon}}-\boldsymbol{\boldsymbol{\varepsilon}}_{s+1}^{p}\right)}{\xi_{s+1}^{\alpha}}\right|^{n}\mathrm{sgn}\left[\tau_{\alpha}\right]=\boldsymbol{0}\label{eq:epdet}
\end{equation}
\\
Newton iteration for $\boldsymbol{\varepsilon}_{s+1}^{p}$ in Voigt form is
performed, which results in a decoupled solution from $\xi_{s+1}^{\alpha}.$
The Jacobian of (\ref{eq:epdet}) is determined as:

\begin{equation}
\boldsymbol{J}_{\varepsilon_{p}}=\boldsymbol{I}_{6\times6}+\Delta t\sum_{\alpha=1}^{12}\boldsymbol{M}_{\alpha}\otimes\boldsymbol{P}_{\alpha}\frac{n\dot{\gamma}_{0}^{\alpha}}{\xi_{s+1}^{\alpha}}\left|\frac{\boldsymbol{P}_{\alpha}:\left(\boldsymbol{\boldsymbol{\varepsilon}}-\boldsymbol{\boldsymbol{\varepsilon}}_{s+1}^{p}\right)}{\xi_{s+1}^{\alpha}}\right|^{n-1}\label{eq:jacob}
\end{equation}

The hardening law is also integrated:

\begin{equation}
\xi_{s+1}^{\alpha}=\xi_{s}^{\alpha}+h_{0}\sum_{\beta=1}^{12}\left[\left|\Delta\gamma_{s+1}^{\beta}\right|\left(1-\frac{\xi_{s+1}^{\beta}}{\xi_{\infty}^{\beta}}\right)h_{\alpha\beta}\right]\label{eq:hard-1}
\end{equation}
Given that (\ref{eq:hard-1}) is a linear system for $\xi_{s+1}^{\gamma}$,
we rewrite it as:

\begin{align}
\xi_{s+1}^{\gamma}\delta_{\alpha\gamma} & =\xi_{s}^{\alpha}+h_{0}\left[\left|\Delta\gamma_{s+1}^{\beta}\right|\left(1-\frac{\delta_{\beta\gamma}\xi_{s+1}^{\gamma}}{\xi_{\infty}^{\beta}}\right)h_{\alpha\beta}\right]\Leftrightarrow\label{eq:hard-1-1}\\
\left(\delta_{\alpha\gamma}+h_{0}h_{\alpha\gamma}\frac{\left|\Delta\gamma_{s+1}^{\gamma}\right|}{\xi_{\infty}^{\gamma}}\right)\xi_{s+1}^{\gamma} & =\xi_{s}^{\alpha}+h_{0}\left|\Delta\gamma_{s+1}^{\beta}\right|h_{\alpha\beta}
\end{align}

If all slipping systems are active, then it follows that $\left|\tau_{\alpha}\right|-\xi_{s+1}^{\alpha}=0$
for $\alpha=1,\ldots,12$. Introducing the sign of $\tau_{\alpha}$
as $s_{\alpha}$, then it follows that $\tau_{\alpha}-s_{\alpha}\xi_{s+1}^{\alpha}=0$
for active systems. For each $\alpha$, a summation in $\beta=1,\ldots,12$
is required:

\begin{equation}
\tau_{\alpha}^{\star}-s_{\alpha}\xi_{s+1}^{\alpha}-\boldsymbol{P}_{\alpha}\cdot\boldsymbol{M}_{\beta}\Delta\gamma_{s+1}^{\beta}=\boldsymbol{0}\label{eq:linearsystem}
\end{equation}

For the sole purpose of calculating $\Delta\gamma_{s+1}^{\alpha}$,
it is convenient to adopt the forward-Euler algorithm for (\ref{eq:hard}),
resulting in a linearized version of $\xi_{s+1}^{\alpha}$, denominated
here as $\xi_{s+1}^{\alpha L}$:

\begin{equation}
\xi_{s+1}^{\alpha L}=\xi_{s}^{\alpha}+T_{\alpha\beta}s_{\beta}\Delta\gamma_{s+1}^{\beta}\label{eq:xisp1}
\end{equation}
\\
where 
\begin{equation}
T_{\alpha\beta}=h_{0}h_{\alpha\beta}\left(1-\frac{\xi_{s}^{\beta}}{\xi_{\infty}^{\beta}}\right)\label{eq:tab}
\end{equation}
and 
\begin{equation}
s_{\beta}=\nicefrac{\tau_{\beta}}{\left|\tau_{\beta}\right|}\label{eq:stau}
\end{equation}

The coupled system, in terms of constitutive unknowns $\{\boldsymbol{\varepsilon}_{s+1}^{p};\xi_{s+1}^{\alpha},\alpha=1,\ldots,12\}$
is written as a system of two equations:

\begin{align}
\boldsymbol{e}_{\varepsilon} & =\boxed{\boldsymbol{\varepsilon}_{s+1}^{p}-\boldsymbol{\varepsilon}_{s}^{p}-\Delta t\sum_{\alpha=1}^{12}\boldsymbol{M}_{\alpha}\dot{\gamma}_{0}^{\alpha}\left|\frac{\boldsymbol{P}_{\alpha}:\left(\boldsymbol{\boldsymbol{\varepsilon}}-\boldsymbol{\boldsymbol{\varepsilon}}_{s+1}^{p}\right)}{\xi_{s+1}^{\alpha}}\right|^{n}\mathrm{sgn}\left[\tau_{\alpha}\right]=\boldsymbol{0}}\label{eq:hard-1-2}\\
e_{\gamma}= & \boxed{\left(\delta_{\alpha\gamma}+h_{0}h_{\alpha\gamma}\frac{\left|\Delta\gamma_{s+1}^{\gamma}\right|}{\xi_{\infty}^{\gamma}}\right)\xi_{s+1}^{\gamma}-\xi_{s}^{\alpha}-h_{0}\left|\Delta\gamma_{s+1}^{\beta}\right|h_{\alpha\beta}=0;\:\gamma=1,\ldots,12}\label{eq:ea6}
\end{align}

Concering the strain measure, we use power-equivalence to obtain a
conjugate stress to the Hencky strain. The second Piola-Kirchhoff
stress $\boldsymbol{S}$ is power-conjugated to $\dot{\boldsymbol{E}}$ and therefore,
in Voigt form,

\begin{align}
\boldsymbol{S}:\dot{\boldsymbol{E}} & =\boldsymbol{\sigma}:\dot{\boldsymbol{\varepsilon}}=\boldsymbol{\sigma}:\frac{\dif\boldsymbol{\boldsymbol{\varepsilon}}}{\dif\boldsymbol{E}}:\dot{\boldsymbol{E}}\Leftrightarrow\nonumber \\
\boldsymbol{S} & =\boldsymbol{\sigma}:\frac{\dif\boldsymbol{\boldsymbol{\varepsilon}}}{\dif\boldsymbol{E}}\label{eq:ssigma}
\end{align}

For moderate strains, the Pad\'e approximation of order $(2,2)$ is
adopted and it is shown to be adequate. The approximation is given
by (see, e.g. \cite{hajidehi2021}):

\begin{equation}
{\displaystyle \log\left[\boldsymbol{X}\right]\cong[2/2]_{\log}(\boldsymbol{\boldsymbol{I}-X})}=3\left(\boldsymbol{X}^{2}-\boldsymbol{I}\right)\left(\boldsymbol{X}^{2}+4\boldsymbol{X}+\boldsymbol{I}\right)^{-1}\label{eq:pade22}
\end{equation}
\\
this ensures the coincidence, at $\boldsymbol{X=I},$ of the approximation
and the function up to the third derivative. Popular alternatives
to (\ref{eq:pade22}) are polynomial approximations with scaling and
squaring algorithm (used for the exponential by Sastre \emph{et al.}
\cite{sastre2015}) and the approximation proposed by Ba\u{z}ant \cite{bazant1998}.
Since $\boldsymbol{\varepsilon}=\frac{1}{2}\log\left[2\boldsymbol{E}+\boldsymbol{I}\right]$,
we have, in matrix form,

\begin{equation}
\boldsymbol{\varepsilon}\cong3\left[\boldsymbol{E}\cdot\boldsymbol{E}+\boldsymbol{E}\right]\cdot\left[2\boldsymbol{E}\cdot\boldsymbol{E}+6\boldsymbol{E}+3\boldsymbol{I}\right]^{-1}\label{eq:approximateepxilon}
\end{equation}

First and second variations of $\boldsymbol{\varepsilon}$ are required. For
the first variation, we have

\begin{align}
\dif\boldsymbol{\varepsilon} & =3\left[\dif\boldsymbol{E}\cdot\boldsymbol{E}+\boldsymbol{E}\cdot\dif\boldsymbol{E}+\dif\boldsymbol{E}\right]\cdot\left[2\boldsymbol{E}\cdot\boldsymbol{E}+6\boldsymbol{E}+3\boldsymbol{I}\right]^{-1}\nonumber \\
 & -3\left[\boldsymbol{E}\cdot\boldsymbol{E}+\boldsymbol{E}\right]\cdot\left[2\boldsymbol{E}\cdot\boldsymbol{E}+6\boldsymbol{E}+3\boldsymbol{I}\right]^{-1}\cdot\left[2\dif\boldsymbol{E}\cdot\boldsymbol{E}+2\boldsymbol{E}\cdot\dif\boldsymbol{E}+6\dif\boldsymbol{E}\right]\cdot\nonumber \\
 & \left[2\boldsymbol{E}\cdot\boldsymbol{E}+6\boldsymbol{E}+3\boldsymbol{I}\right]^{-1}\label{eq:difepsilon}
\end{align}

The material tangent modulus, $\mathbb{C}$, is calculated using the
chain rule. Using Voigt form,

\begin{equation}
\mathbb{C}_{\textrm{Voigt}}=\left[\frac{\dif\boldsymbol{\boldsymbol{\varepsilon}}}{\dif\boldsymbol{E}}\right]^{T}\cdot\frac{\dif\boldsymbol{\sigma}}{\dif\boldsymbol{\varepsilon}}\cdot\frac{\dif\boldsymbol{\varepsilon}}{\dif\boldsymbol{E}}+\boldsymbol{\sigma}\cdot\frac{\dif^{2}\boldsymbol{\boldsymbol{\varepsilon}}}{\dif\boldsymbol{E}\dif\boldsymbol{E}}\label{eq:tangmodulus}
\end{equation}

The specific expression for $\mathbb{C}_{\textrm{Voigt}}$ is too
intricate to present in this text. The corresponding Mathematica/Acegen
source code is available in the Github repository \cite{areias2023github}.
The optimized expression for (\ref{eq:tangmodulus}) is compact (around
840 lines of Fortran 95) and dispenses the explicit calculation of
the sixth-order tensor $\nicefrac{\dif^{2}\boldsymbol{\varepsilon}}{\dif\boldsymbol{E}\dif\boldsymbol{E}}$
described in Miehe, Apel and Lambrecht \cite{miehe2002}. Error analysis
by C. Kenney and A. Laub \cite{kenney1989} presents the following
inequality, relating the errors of Pad\'e approximation (after specialization
for the present case):

\begin{equation}
\left\Vert [m/n]_{\log}(-2\boldsymbol{\boldsymbol{E}})-\log\left[2\boldsymbol{E}+\boldsymbol{I}\right]\right\Vert \leq\underbrace{[m/n]_{\log}(2\boldsymbol{\left\Vert \boldsymbol{\boldsymbol{E}}\right\Vert })-\log\left[1-2\left\Vert \boldsymbol{\boldsymbol{E}}\right\Vert \right]}_{2E_{[m/n]\log}^{\max}}\label{eq:ineq}
\end{equation}

\begin{figure}
\begin{centering}
\includegraphics[clip,width=11cm]{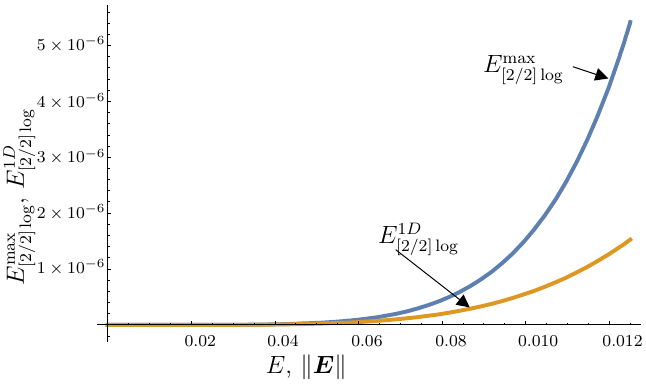}
\par\end{centering}
\caption{\label{figastaggered}Upper bound on the error $E_{[2/2]\log}^{\max}$
compared with the error in closed form $E_{[2/2]\log}^{1D}$.}
\end{figure}

This provides an upper bound for the error in $\boldsymbol{\varepsilon}$.
Bounds for the condition number are provided in \cite{kenney1989}.
Of course, in the 1D case, the absolute error can be calculated in
closed form:

\begin{equation}
E_{[m,n]\log}^{1D}=\nicefrac{1}{2}\left|[m/n]_{\log}(2\boldsymbol{E})-\log\left[1+2\boldsymbol{E}\right]\right|\label{eq:oned}
\end{equation}
\\
A graphical representation of $E_{[2/2]\log}^{\max}$ as a function
of $\|\boldsymbol{E}\|$ and of $E_{[2,2]\log}^{1D}$, respectively, is shown
in Figure \ref{figastaggered}. For materials with limited elastic
strains, the error in the approximation of logarithm is compatible
with current computational mechanics practice.

For the incompressible case, which is reasonable in metal plasticity
with finite strains, and using principal directions (principal Euler-Lagrange
strains $E_{1}$, $E_{2}$ and $E_{3}$), we have, for the Hencky
strain,

\begin{equation}
\nicefrac{1}{2}\log\left[-2\boldsymbol{E}\right]=\nicefrac{1}{2}\left[\begin{array}{ccc}
\log\left[2E_{1}+1\right] & 0 & 0\\
0 & \log\left[2E_{2}+1\right] & 0\\
0 & 0 & -\log\left[1+2E_{1}+2E_{2}+4E_{1}E_{2}\right]
\end{array}\right]\label{eq:logstrain}
\end{equation}
\\
where it was assumed that $E_{3}$ is a function of $E_{1}$ and $E_{2}$:

\begin{equation}
E_{3}=\frac{1}{2}\left[-1+\frac{1}{1+2E_{1}+2E_{2}+4E_{1}E_{2}}\right]\label{eq:e3det}
\end{equation}
For the $[2/2]$ Pad\'e approximation, 

\begin{equation}
\boldsymbol{\varepsilon}=\left[\begin{array}{ccc}
\frac{3E_{1}\left(1+E_{1}\right)}{3+2E_{1}\left(3+E_{1}\right)} & 0 & 0\\
0 & \frac{3E_{2}\left(1+E_{2}\right)}{3+2E_{2}\left(3+E_{2}\right)} & 0\\
0 & 0 & \frac{3E_{3}\left(1+E_{3}\right)}{3+2E_{3}\left(3+E_{3}\right)}
\end{array}\right]\label{eq:epsilon}
\end{equation}

\begin{figure}
\begin{centering}
\includegraphics[width=8cm]{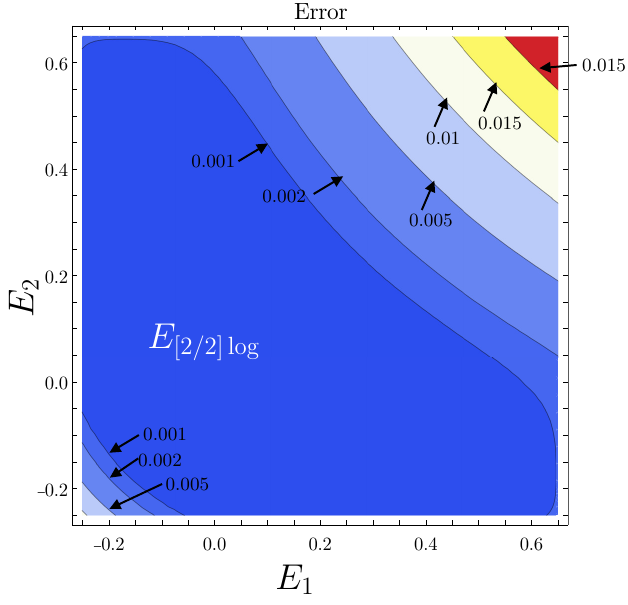}
\par\end{centering}
\caption{\label{figaRelative-error-for}Relative error for an incompressible
$2D$ problem, $E_{[2,2]\log}=\left\Vert \nicefrac{1}{2}[2/2]_{\log}(-2\boldsymbol{E})-\boldsymbol{\varepsilon}\right\Vert $.}

\end{figure}

Figure \ref{figaRelative-error-for} shows the error $\|\nicefrac{1}{2}\log\left[2\boldsymbol{E}+I\right]-\boldsymbol{\varepsilon}\|$
in the domain $E_{k}\in[-0.25,0.65]$. It can be observed that even
for considerable strains in tension, the error is $0.0258$, which
corresponds to a relative error of $2.52$\%.

The calculation of the rotated normals is performed with Kr\"{o}ner/Lee's
multiplicative decomposition \cite{kroner1960,lee1967,lee1969} of
the deformation gradient $\boldsymbol{F}$ into elastic $\boldsymbol{F}^{e}$ and
plastic parts $\boldsymbol{F}^{p}$, 

\begin{equation}
\boldsymbol{F}=\boldsymbol{F}^{e}\cdot\boldsymbol{F}^{p}\label{eq:fefp}
\end{equation}
\\
Since $\boldsymbol{F}$ is determined by the finite element solution, $\boldsymbol{F}_{e}$
at time step $s$ is determined by the flow law (see, e.g. \cite{alankar2009})
as:

\begin{equation}
\boldsymbol{F}_{s}^{e}=\boldsymbol{F}_{s}\cdot\left[\boldsymbol{F}_{s-1}^{p}\right]^{-1}\cdot\left(\boldsymbol{I}-\sum_{\alpha=1}^{12}\boldsymbol{m}_{\alpha}\otimes\boldsymbol{n}_{\alpha}\Delta t\dot{\gamma}^{\alpha}\right)\label{eq:fe}
\end{equation}
\\
from which the plane normals in the deformed configuration (see Kaiser
and Menzel \cite{kaiser2019}) are obtained using the elastic push-forward:

\begin{equation}
\boldsymbol{m}_{\alpha}^{\star}=\boldsymbol{F}_{s}^{e}\cdot\boldsymbol{m}_{\alpha}\label{eq:fepullback}
\end{equation}

\begin{figure}
\begin{centering}
\subfloat[Geometry and boundary conditions. Time increment is $\Delta t=1\times10^{-3}$
s and $\dot{\overline{u}}_{3}=0.001$ $\textrm{m/s}$]{\begin{centering}
\includegraphics[width=9cm]{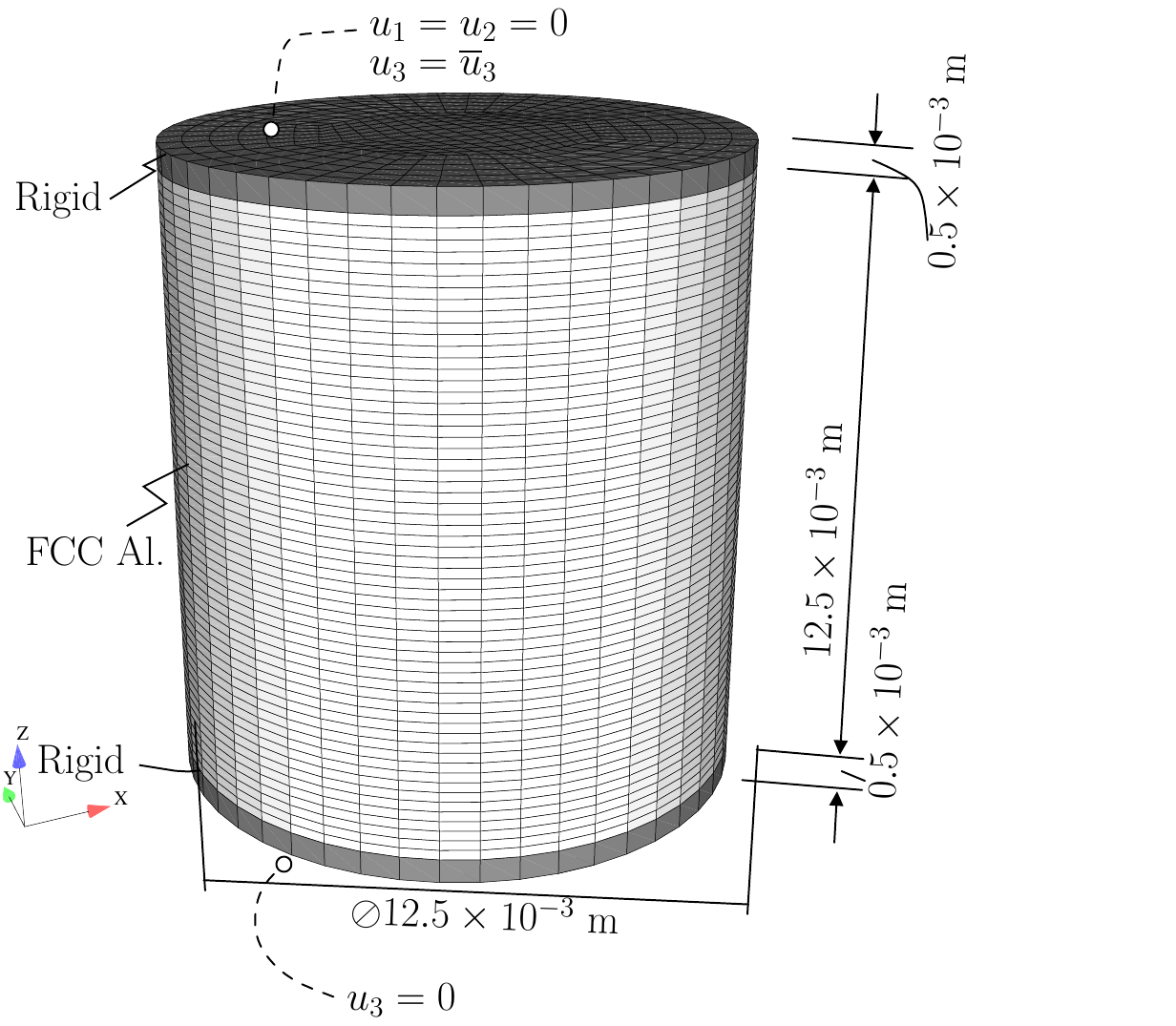}
\par\end{centering}
}
\par\end{centering}
\begin{centering}
\subfloat[Cases of plane $(111)$ orientation with respect to $z$]{\begin{centering}
\includegraphics[width=9cm]{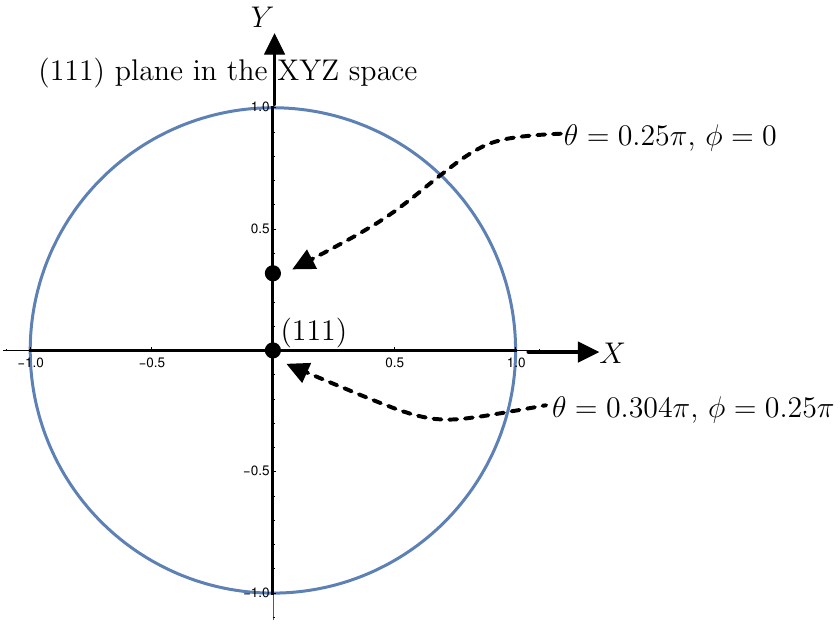}
\par\end{centering}

}
\par\end{centering}
\caption{\label{figadefmmmeshes-1}Verification test: geometry and boundary
conditions for the single crystal cylinder under tension. }
\end{figure}

We use the centroidal deformation gradient for purposes of calculating
$\boldsymbol{F}_{s}^{e}$ and $\boldsymbol{m}_{\alpha}^{\star}$. A verification test
is performed using the data shown in Figure \ref{figadefmmmeshes-1}.
Two rigid plates are connected to the single crystal. The upper plate
is clamped and pulled in the $z$ direction by an imposed displacement
$\overline{u}_{3}$ where $\dot{\overline{u}}_{3}=0.001$ $\textrm{m/s}$.
The lower plate is fixed in the $z$ direction but left free to have
displacement in the $x-y$ plane. We measure the average normal stress
$\overline{\sigma}_{33}$ by dividing the reaction force on the $z$
direction. In addition, the average strain, $\overline{E}_{33}$ is
obtained as $\overline{E}_{33}=\nicefrac{\overline{u}_{3}}{12.5\times10^{-3}}.$
To force strain localization, we adopt the saturation stress $\xi_{\infty}^{\star}$
in Table \ref{tab:Hardening-properties}. 

\begin{figure}
\begin{centering}
\includegraphics[width=14cm]{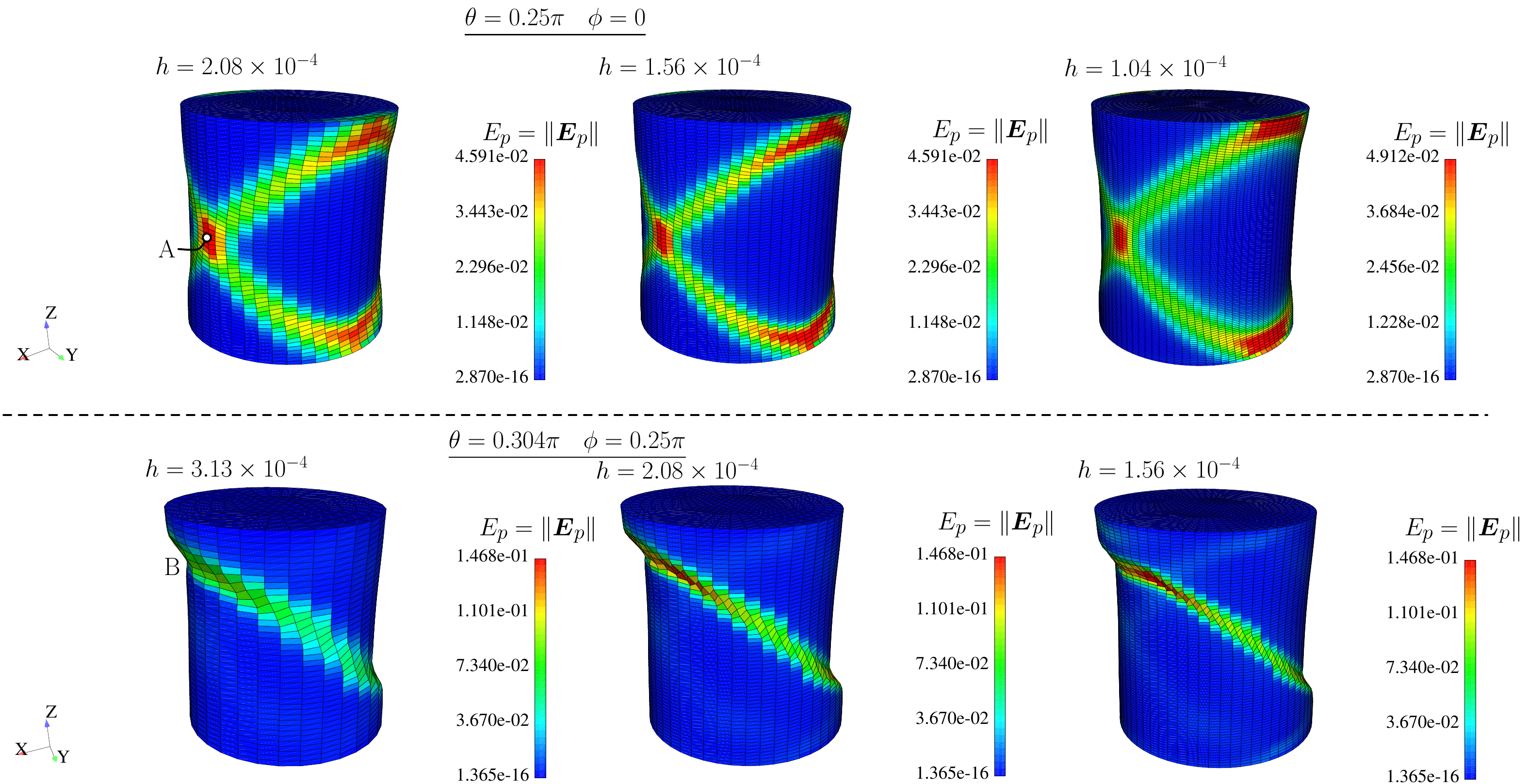}
\par\end{centering}
\caption{\label{figadefmmmeshes}Deformed meshes ($10\times$ magnified)
and $E_{p}=\|\boldsymbol{E}_{p}\|$ contour plots for $\{\theta\}$. Three
characteristic mesh sizes $h=3.13\times10^{-4},2.08\times10^{-4},1.56\times10^{-4}$
are tested for $\theta=0.304\pi$ and $\phi=0.25\pi$ and $h=2.08\times10^{-4},1.56\times10^{-4},1.04\times10^{-4}$
for $\theta=0.25\pi$ and $\phi=0$. Consistent units are used.}
\end{figure}

\begin{figure}
\begin{centering}
\subfloat[$\theta=0.25\pi$, $\phi=0$]{\begin{centering}
\includegraphics[width=11cm]{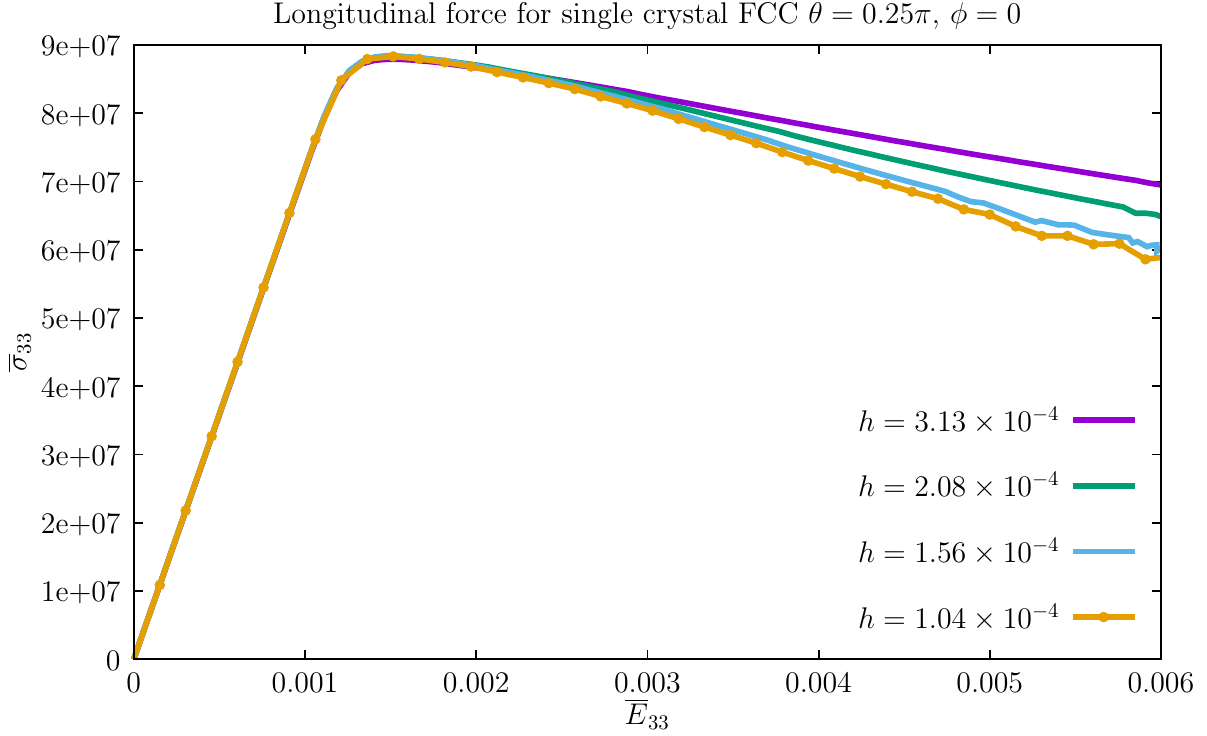}
\par\end{centering}
}
\par\end{centering}
\begin{centering}
\subfloat[$\theta=0.304\pi$, $\phi=0.25\pi$]{\begin{centering}
\includegraphics[width=11cm]{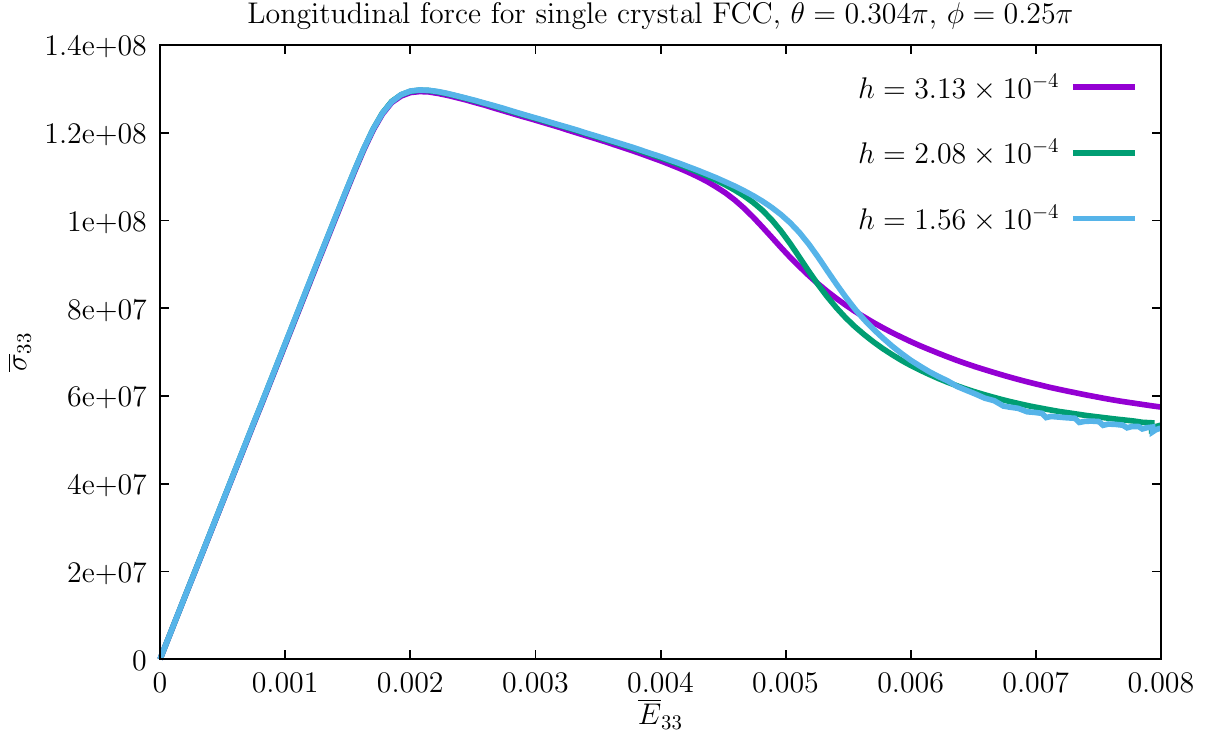}
\par\end{centering}
}
\par\end{centering}
\caption{\label{figaEffect-of-}Effect of $h$ on the average stress $\overline{\sigma}_{33}$
and $\overline{E}_{33}$ for two orientations: $\theta=0.25\pi$,
$\phi=0$ and $\theta=0.304\pi$, $\phi=0.25\pi$.}

\end{figure}

Since strain rate dependence is present in the flow law, strain softening
is allowed. Localization results for the two orientations are shown
in Figure \ref{figadefmmmeshes} as a function of the element
size $h$ (see, e.g. \cite{hughes2000} for this nomenclature). The
relation between these quantities is represented in Figure \ref{figaEffect-of-},
where good mesh insensitivity can be observed, especially considering
that no regularization is adopted. For $\theta=0.304\pi$ and $\phi=0.25\pi$,
we show the contour plots of $\xi_{1},\cdots,\xi_{12}$ in Figure
\ref{figaContour-plots-of}.

\begin{center}
\begin{figure}
\begin{centering}
\includegraphics[width=14cm]{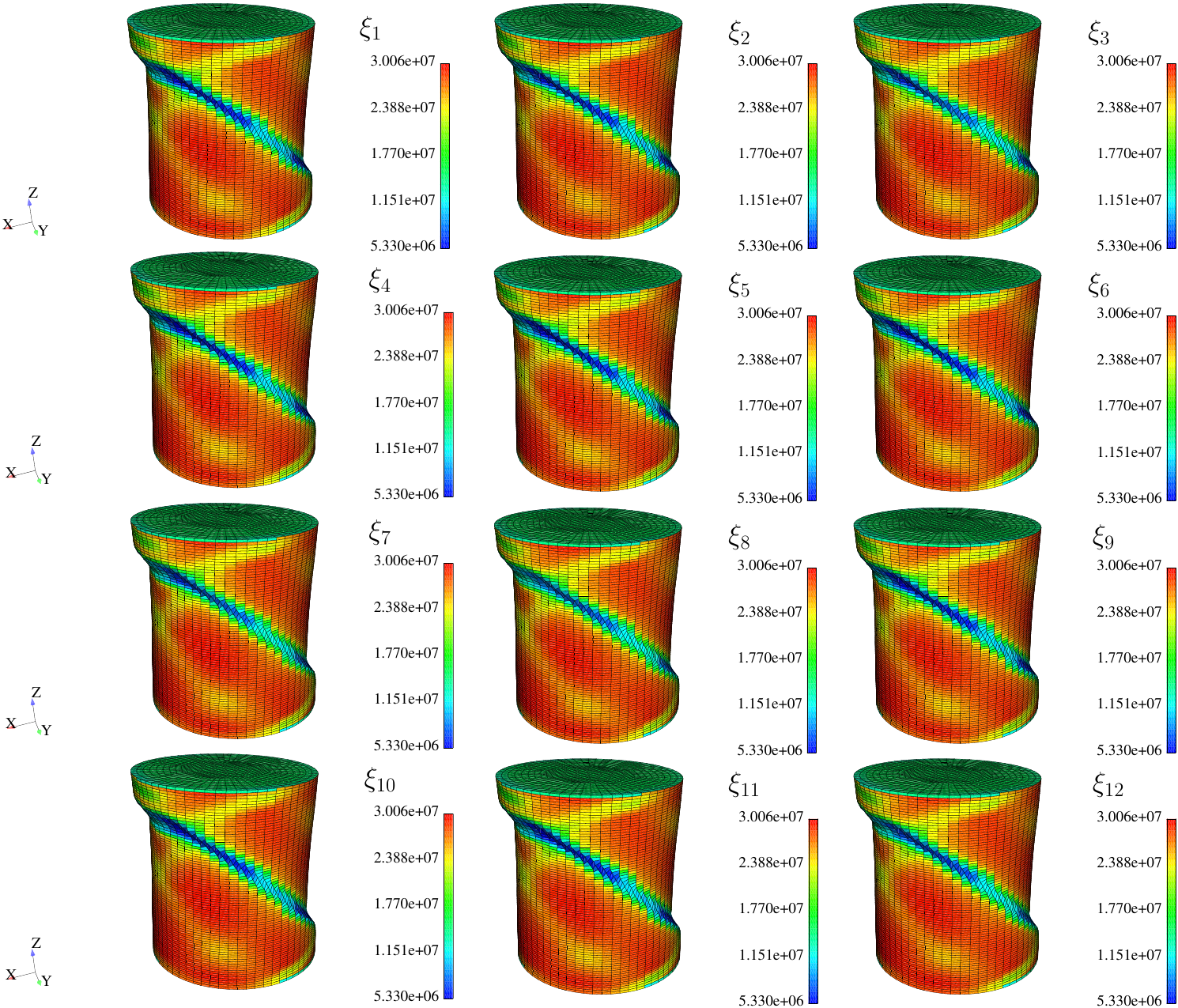}
\par\end{centering}
\caption{\label{figaContour-plots-of}Contour plots of $\xi_{1},\cdots,\xi_{12}$
for $\theta=0.304\pi$,  $\phi=0.25\pi$.}

\end{figure}
\par\end{center}
\section{Finite element formulation}

\label{secfiniteelementformulation}

We adopt the 8-node hexahedron with EAS (enhanced assumed strain)
technology developed by Simo and co-workers \cite{simorifai1990,simoarmero1992,simo1993},
which is appropriate for strain localization problems (cf. \cite{simoarmero1992}).
A version based on the Green-Lagrange strain was introduced by Andelfinger
and Ramm \cite{andelfinger1993} and it is employed here. Calculations
were performed using Mathematica with the AceGen add-on \cite{korelc2002}.
Using the classical formalism and the Euler-Lagrange strains, 

\begin{equation}
\boldsymbol{E}=\boldsymbol{E}_{u}+\boldsymbol{E}_{\alpha}\label{eq:edecomp}
\end{equation}
\\
see, e.g. \cite{andelfinger1993}. Classical results for $\boldsymbol{E}_{u}$
and $\boldsymbol{E}_{\alpha}$ follow (see, e.g. \cite{andelfinger1993}):

\begin{equation}
\boldsymbol{E}_{u}=\frac{1}{2}\left[\left(\nabla\boldsymbol{u}\right)+\left(\nabla\boldsymbol{u}\right)^{T}+\left(\nabla\boldsymbol{u}\right)^{T}\cdot\left(\nabla\boldsymbol{u}\right)\right]\label{eq:eu}
\end{equation}

\begin{equation}
\boldsymbol{E}_{\alpha}=\frac{\overline{J}}{J}\overline{\boldsymbol{J}}^{-T}\cdot\widetilde{\boldsymbol{E}}_{\alpha}\cdot\overline{\boldsymbol{J}}^{-1}\label{eq:eeas}
\end{equation}
\\
where $\left[\nabla\boldsymbol{u}\right]_{ij}=\left[\nicefrac{\partial u_{i}}{\partial X_{j}}\right]$
is the displacement gradient with respect to the undeformed coordinates
and $\widetilde{\boldsymbol{E}}_{\alpha}$ is the enhanced strain in parent-domain
(see, e.g. \cite{simo1993}) coordinates. It is required to adopt
the centroidal Jacobian matrix $\overline{\boldsymbol{J}}$ and the corresponding
determinant $\overline{J}$ to ensure the patch-test satisfaction.
The EAS modes are based on the interpolation of the bending deformation
modes. $12$ additional, internal element degrees-of-freedom are adopted,
which here we denote by the letter $\boldsymbol{\alpha}=\left\{ \alpha_{1},\ldots,\alpha_{12}\right\} $.
Using Voigt form \cite{belytschko2000}, it follows that:

\[
\widetilde{\boldsymbol{E}}_{\alpha}^{\textrm{voigt}}=\left[\begin{array}{c|c|c|c|c|c|c|c|c|c|c|c}
\xi_{2} & \xi_{3} & \xi_{2}\xi_{3} & 0 & 0 & 0 & 0 & 0 & 0 & 0 & 0 & 0\\
\hline 0 & 0 & 0 & \xi_{1} & \xi_{3} & \xi_{1}\xi_{3} & 0 & 0 & 0 & 0 & 0 & 0\\
\hline 0 & 0 & 0 & 0 & 0 & 0 & \xi_{1} & \xi_{2} & \xi_{1}\xi_{2} & 0 & 0 & 0\\
\hline 0 & 0 & 0 & 0 & 0 & 0 & 0 & 0 & 0 & \xi_{3} & 0 & 0\\
\hline 0 & 0 & 0 & 0 & 0 & 0 & 0 & 0 & 0 & 0 & \xi_{2} & 0\\
\hline 0 & 0 & 0 & 0 & 0 & 0 & 0 & 0 & 0 & 0 & 0 & \xi_{1}
\end{array}\right]\cdot\left\{ \begin{array}{c}
\alpha_{1}\\
\hline \vdots\\
\hline \alpha_{12}
\end{array}\right\} 
\]

As mentioned, patch test is satisfied a-priori by the the use of $\overline{J}$
and $\overline{\boldsymbol{J}}$ in (\ref{eq:eeas}), see also \cite{simoarmero1992,simo1993}.
The full weak form corresponding to the decomposition is as follows,
using the undeformed configuration:

\begin{equation}
\underbrace{\int_{\Omega_{0}}\boldsymbol{S}:\left(\delta\boldsymbol{E}_{u}+\delta\boldsymbol{E}_{\alpha}\right)\dif V}_{\delta W_{\textrm{int}}}=\underbrace{\int_{\Omega_{0}}\rho_{0}\overline{B}\,\cdot\delta\boldsymbol{u}\dif V+\int_{\partial\Omega_{0}}\overline{\boldsymbol{T}}\cdot\delta\boldsymbol{u}\dif A}_{\delta W_{\textrm{ext}}}\label{eq:se}
\end{equation}
\\
where the Piola stress vector $\overline{\boldsymbol{T}}$, resulting from
integration by parts, is calculated as (see \cite{wriggers2008}):
$\overline{\boldsymbol{T}}=\boldsymbol{F}\cdot\boldsymbol{S}\cdot\boldsymbol{N}$ where $\boldsymbol{F}$
is the deformation gradient and $\boldsymbol{N}$ is the outer normal to $\partial\Omega_{0}.$
The second Piola-Kirchhoff stress $\boldsymbol{S}$ is a function of the total
Green-Lagrange strain $\boldsymbol{E}$ as $\boldsymbol{S}\equiv\boldsymbol{S}\left(\boldsymbol{E}_{u}+\boldsymbol{E}_{\alpha}\right)$.
Use of Newton iteration for (\ref{eq:se}) requires the calculation
of the variation in both hand-sides:

\begin{equation}
\Delta_{u}\delta W_{\textrm{int}}\cdot\Delta\boldsymbol{u}+\Delta_{\alpha}\delta W_{\textrm{int}}\cdot\Delta\boldsymbol{\alpha}=\delta W_{\textrm{ext}}-\delta W_{\textrm{int}}\label{eq:wei}
\end{equation}
\\
assuming that $\Delta_{u}\delta W_{\textrm{ext }}=0$. Source code
for this element (forces and tangent stiffness) is available in Github
\cite{areias2024element}. It is worth noting that a development of
EAS has been used with single crystal plasticity by J. Mosler's group
\cite{fohrmeister2019}.

\section{Implicit staggered algorithm}

\label{secimplstagg}

Since a decoupled constitutive system (\ref{eq:hard-1-2}-\ref{eq:ea6})
is solved, corresponding to the flow law and the hardening evolution,
explicit algorithms produce drift, as is the case in thermoelasticity
\cite{erbts2012}, fluid-structure interaction (FSI) \cite{kuttler2008,degroote2010}
and phase-field simulations \cite{schapira2023}. The study by Erbts
and D\"{u}ster shows that dynamic relaxation with appropriate predictors
produces efficient and stable results to remove the drift. We adopt
dynamic relaxation and introduce a substep index $i$ as a superscript.
Hardening variables $\boldsymbol{\xi}$ are adopted in the relaxation. Index
interpretation is as follows:

\[
\boldsymbol{\xi}_{s}^{i}:\quad i^{\textrm{th}}\textrm{\,substep of}\quad\textrm{time\ensuremath{\quad\textrm{step \ensuremath{s}}}}
\]

Since the value of $\boldsymbol{\xi}$ depends on $\Delta\boldsymbol{\gamma}$ and
this also depends on $\boldsymbol{\xi}$, the staggered solution must comply
with the original coupled system. To represent these dependencies,
we introduce the operator $\Xi_{\star}(\boldsymbol{\xi})$ as:

\begin{equation}
\widetilde{\boldsymbol{\xi}}_{s+1}^{i+1}=\Xi_{\star}\left(\boldsymbol{\xi}_{s+1}^{i}\right)=\Xi\Bigl[\underbrace{\Gamma\left(\boldsymbol{\xi}_{s+1}^{i}\right)}_{\Delta\boldsymbol{\gamma}},\boldsymbol{\xi}_{s+1}^{i}\Bigr]\label{eq:xitil}
\end{equation}
\\
,~see also Eq. (28) in \cite{erbts2012} for a similar operator.
Relaxation methods make use of a combination of fixed-point iteration
with heuristic acceleration. Updating is based on linear combination
of the image (\ref{eq:xitil}) and the previous substep:

\begin{equation}
\boldsymbol{\xi}_{s+1}^{i+1}=\left(1-\omega_{i}\right)\boldsymbol{\xi}_{s+1}^{i}+\omega_{i}\widetilde{\boldsymbol{\xi}}_{s+1}^{i+1}\label{eq:xisi}
\end{equation}
\\
where $\omega_{i}$ is the coefficient of the linear combination,
$\omega_{i}\in]0,2[$. The hardening variable residual for substep
$i+1$ is given by:

\begin{equation}
\boldsymbol{r}_{i}=\widetilde{\boldsymbol{\xi}}_{s+1}^{i+1}-\boldsymbol{\xi}_{s+1}^{i}\label{eq:ri}
\end{equation}

Heuristics for updating $\omega_{i}$ have been discussed at length.
We here follow \cite{erbts2012} to update $\omega_{i}:$

\begin{equation}
\omega_{i}=\omega_{i-1}\left[1+\frac{\left(\boldsymbol{r}_{i-1}-\boldsymbol{r}_{i}\right)\cdot\boldsymbol{r}_{i}}{\left(\boldsymbol{r}_{i-1}-\boldsymbol{r}_{i}\right)\cdot\left(\boldsymbol{r}_{i-1}-\boldsymbol{r}_{i}\right)}\right]\label{eq:wii}
\end{equation}

Figure \ref{figastaggered} shows the algorithm \ref{alg:Anderson-algorithm-for}
adopted in our code.
\begin{center}
\begin{algorithm}
\caption{\label{alg:Anderson-algorithm-for}Dynamic relaxation algorithm for
$\Delta\boldsymbol{\gamma}$ and $\boldsymbol{\xi}$}

\vspace{0.3cm}
\begin{centering}
\begin{tabular}{c|c|ccc}
\cline{1-3} \cline{2-3} \cline{3-3} 
$1$ &  & $\boldsymbol{\xi}_{s+1}^{0}=\boldsymbol{\xi}_{s}$ &  & \tabularnewline
$2$ &  & Newton iteration for $\Delta\boldsymbol{\gamma}$: $\Delta\boldsymbol{\gamma}=\Gamma\left(\boldsymbol{\xi}_{s+1}^{0}\right)$ &  & \tabularnewline
$3$ &  & Linear solution for $\widetilde{\boldsymbol{\xi}}_{s+1}^{1}=\Xi\left[\Delta\boldsymbol{\gamma},\boldsymbol{\xi}_{s+1}^{0}\right]$ &  & \tabularnewline
$4$ &  & Initial substep residual $\boldsymbol{r}_{0}=\widetilde{\boldsymbol{\xi}}_{s+1}^{1}-\boldsymbol{\xi}_{s+1}^{0}$ &  & \tabularnewline
$5$ &  & Initial coefficient $\omega_{0}=\nicefrac{1}{2}$ &  & \tabularnewline
$6$ &  & Estimate $\boldsymbol{\xi}_{s+1}^{1}=\left(1-\omega_{0}\right)\boldsymbol{\xi}_{s+1}^{0}+\omega_{0}\widetilde{\boldsymbol{\xi}}_{s+1}^{1}$ &  & \tabularnewline
$7$ &  & Iteration $i=1,\ldots$until $\|\boldsymbol{r}_{i+1}\|\leq\varepsilon\|\boldsymbol{r}_{0}\|$ &  & \tabularnewline
 & $7.1$ & $\boldsymbol{\widetilde{\xi}}_{s+1}^{i+1}\leftarrow\Xi_{\star}\left(\boldsymbol{\xi}_{s+1}^{i}\right)$ &  & \tabularnewline
 & $7.2$ & $\boldsymbol{r}_{i}=\widetilde{\boldsymbol{\xi}}_{s+1}^{i+1}-\boldsymbol{\xi}_{s+1}^{i}$ &  & \tabularnewline
 & $7.3$ & $\|\boldsymbol{r}_{i}\|\leq\varepsilon\|\boldsymbol{r}_{0}\|\rightarrow\textrm{goto 8}$ &  & \tabularnewline
 & $7.4$ & $\omega_{i}=\omega_{i-1}\left[1+\frac{\left(\boldsymbol{r}_{i-1}-\boldsymbol{r}_{i}\right)\cdot\boldsymbol{r}_{i}}{\left(\boldsymbol{r}_{i-1}-\boldsymbol{r}_{i}\right)\cdot\left(\boldsymbol{r}_{i-1}-\boldsymbol{r}_{i}\right)}\right]$ &  & \tabularnewline
 & $7.5$ & $\boldsymbol{\xi}_{s+1}^{i+1}=\left(1-\omega_{i}\right)\boldsymbol{\xi}_{s+1}^{i}+\omega_{i}\widetilde{\boldsymbol{\xi}}_{s+1}^{i+1}$ &  & \tabularnewline
 & $7.6$ & $i\leftarrow i+1$ &  & \tabularnewline
$8$ &  & End iteration &  & \tabularnewline
\cline{1-3} \cline{2-3} \cline{3-3} 
\multicolumn{1}{c}{} & \multicolumn{1}{c}{} &  &  & \tabularnewline
\end{tabular}
\par\end{centering}
\vspace{0.3cm}
\end{algorithm}
\par\end{center}

The effect of time step size $\Delta t$ in the drifting is assessed
in Figure \ref{figaEffect-of-dynamic}. The uncoupled version is run
with two successive passes for $\Delta\boldsymbol{\gamma}$ and $\boldsymbol{\xi}_{s+1}$
and this is compared with Algorithm \ref{alg:Anderson-algorithm-for}.
We can 
\begin{center}
\begin{figure}
\begin{centering}
\subfloat[{Evolution of $\left[\varepsilon_{p}\right]_{33}$ for $\theta=0.25\pi$,
$\phi=0$ in point A.}]{\begin{centering}
\includegraphics[clip,width=10cm]{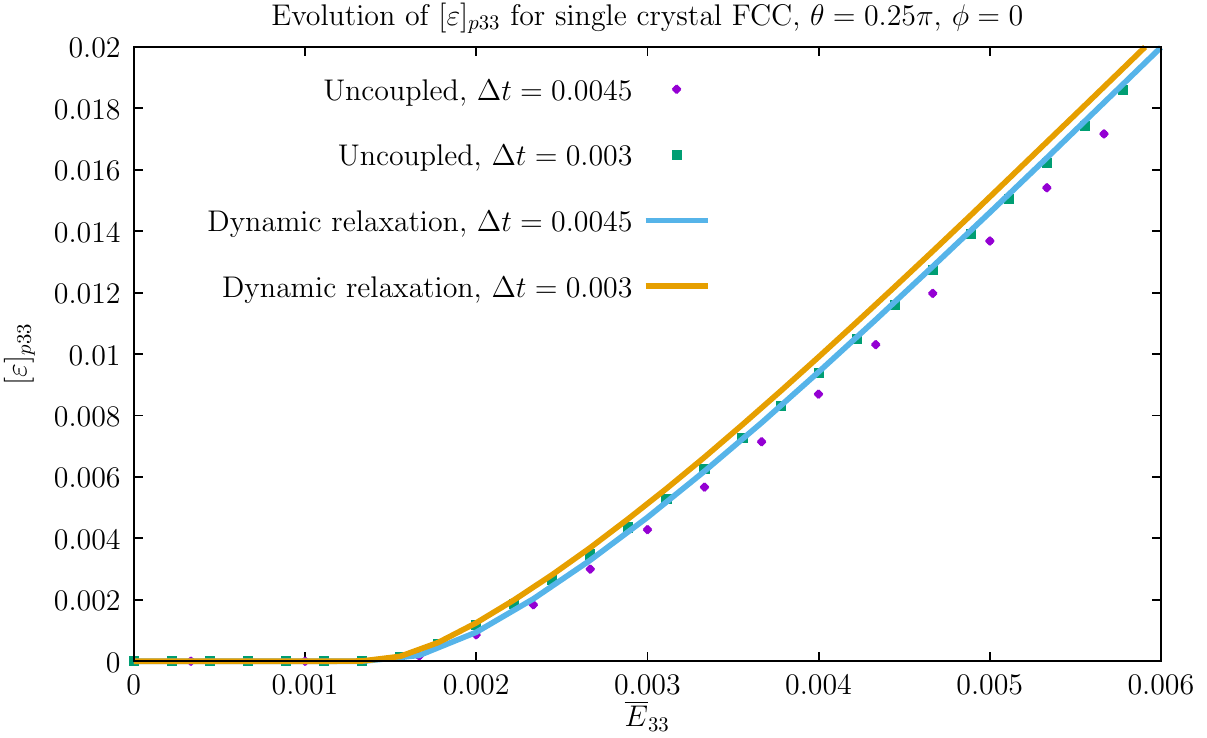}
\par\end{centering}
}
\par\end{centering}
\begin{centering}
\subfloat[{Evolution of $\left[\varepsilon_{p}\right]_{33}$ for $\theta=0.304\pi$,
$\phi=0.25\pi$ in point B.}]{\begin{centering}
\includegraphics[clip,width=10cm]{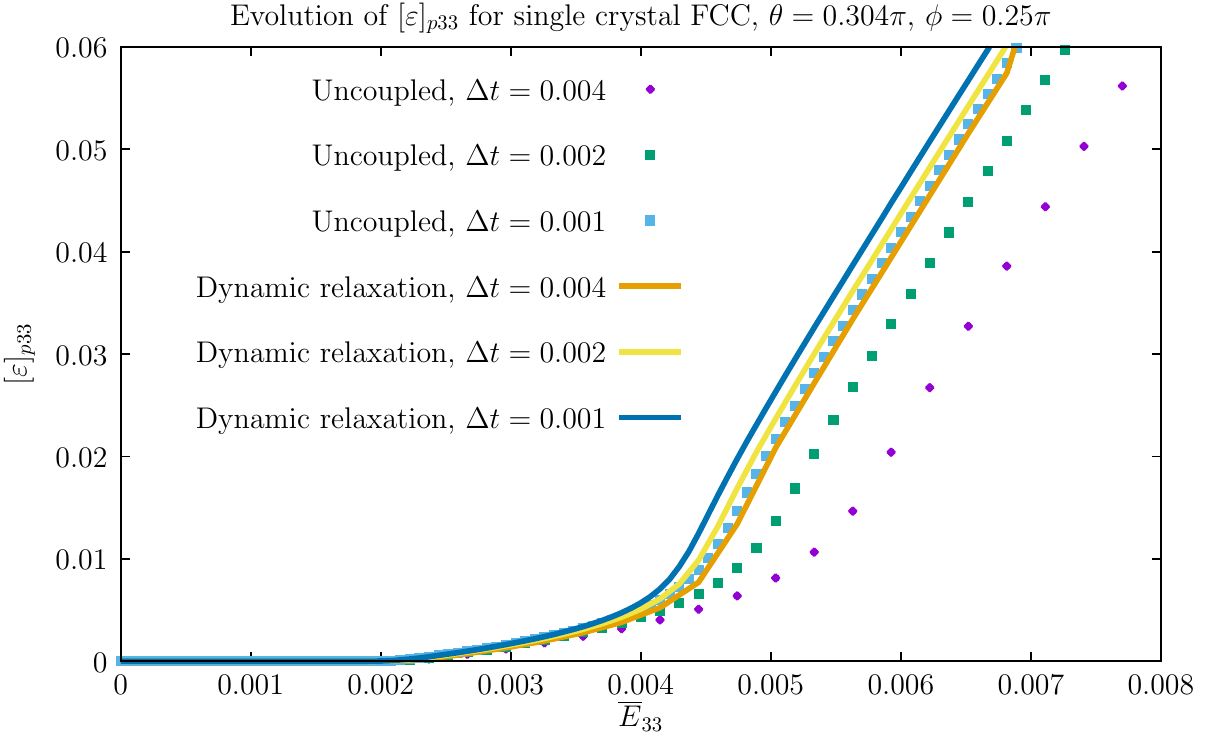}
\par\end{centering}
 }
\par\end{centering}
\caption{\label{figaEffect-of-dynamic}Effect of dynamic relaxation on the
evolution of $\left[\varepsilon_{p}\right]_{33}$ .}
\end{figure}
\par\end{center}

The algorithm is inserted in our in-house code, SimPlas \cite{simplascode}
and combined with the existing finite element technology.

\section{Numerical assessment}

\label{secnumass}

\begin{figure}
\begin{centering}
\subfloat[Geometry and boundary conditions]{\begin{centering}
\includegraphics[width=12cm]{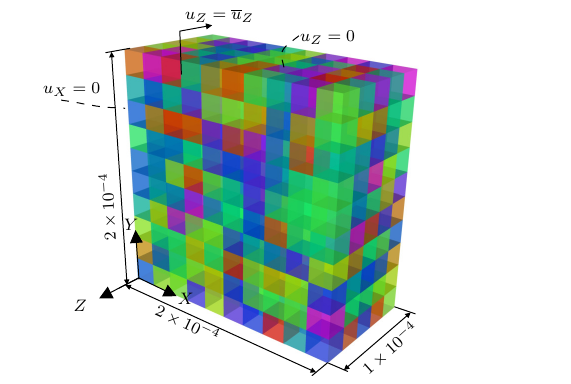}
\par\end{centering}
}
\par\end{centering}
\begin{centering}
\subfloat[\label{figaDistribution-of-(111)}Distribution of (111) planes in
the undeformed configuration.]{\begin{centering}
\includegraphics[width=7cm]{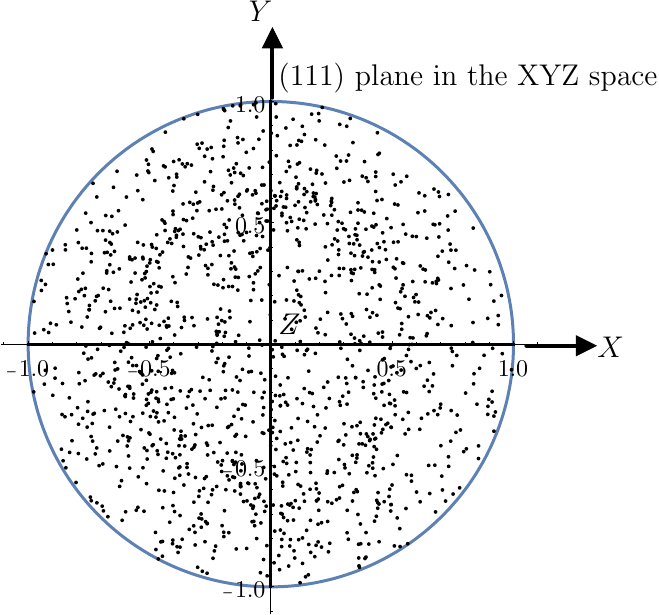}
\par\end{centering}

}
\par\end{centering}
\caption{\label{figaAluminum-polycrystal}Aluminum polycrystal: geometry, dimensions
and boundary conditions. Also shown is the distribution of $(111)$
planes in the undeformed configuration.}

\end{figure}

\begin{figure}
\begin{centering}
\includegraphics[width=14cm]{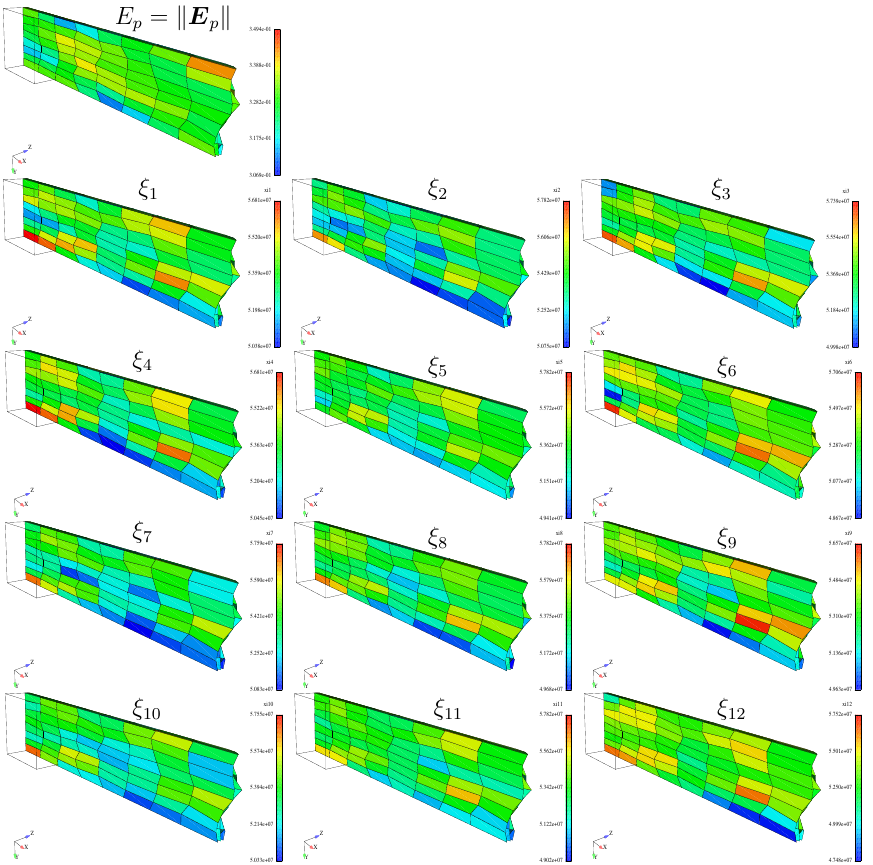}
\par\end{centering}
\caption{\label{figaAluminum-polycrystal:-distributi}Aluminum polycrystal:
distribution of $E_{p}=\|\boldsymbol{E}_{p}\|$ and $\xi_{1},\cdots,\xi_{12}$
for $80\%$ compression.}

\end{figure}

We now test a polycrystalline rectangular cuboid as depicted in Figure
\ref{figaAluminum-polycrystal}. The initial distribution of $(111)$
plane in the XYZ space is shown in sub-Figure \ref{figaDistribution-of-(111)}.
This test was proposed by Alankar, Mastorakos and Field \cite{alankar2009}.
In contrast with that reference, a free edge with $X=2\times10^{-4}$
is left to exhibit the texture. The contour plots of $\xi_{k}$ and
$E_{p}=\|\boldsymbol{E}_{p}\|$ for $80\%$ compression are shown in Figure
\ref{figaAluminum-polycrystal:-distributi}. These are close to those
reported in \cite{alankar2009} but the right edge $X=2\times10^{-4}$
is here left free.

\begin{figure}
\begin{centering}
\includegraphics[width=12cm]{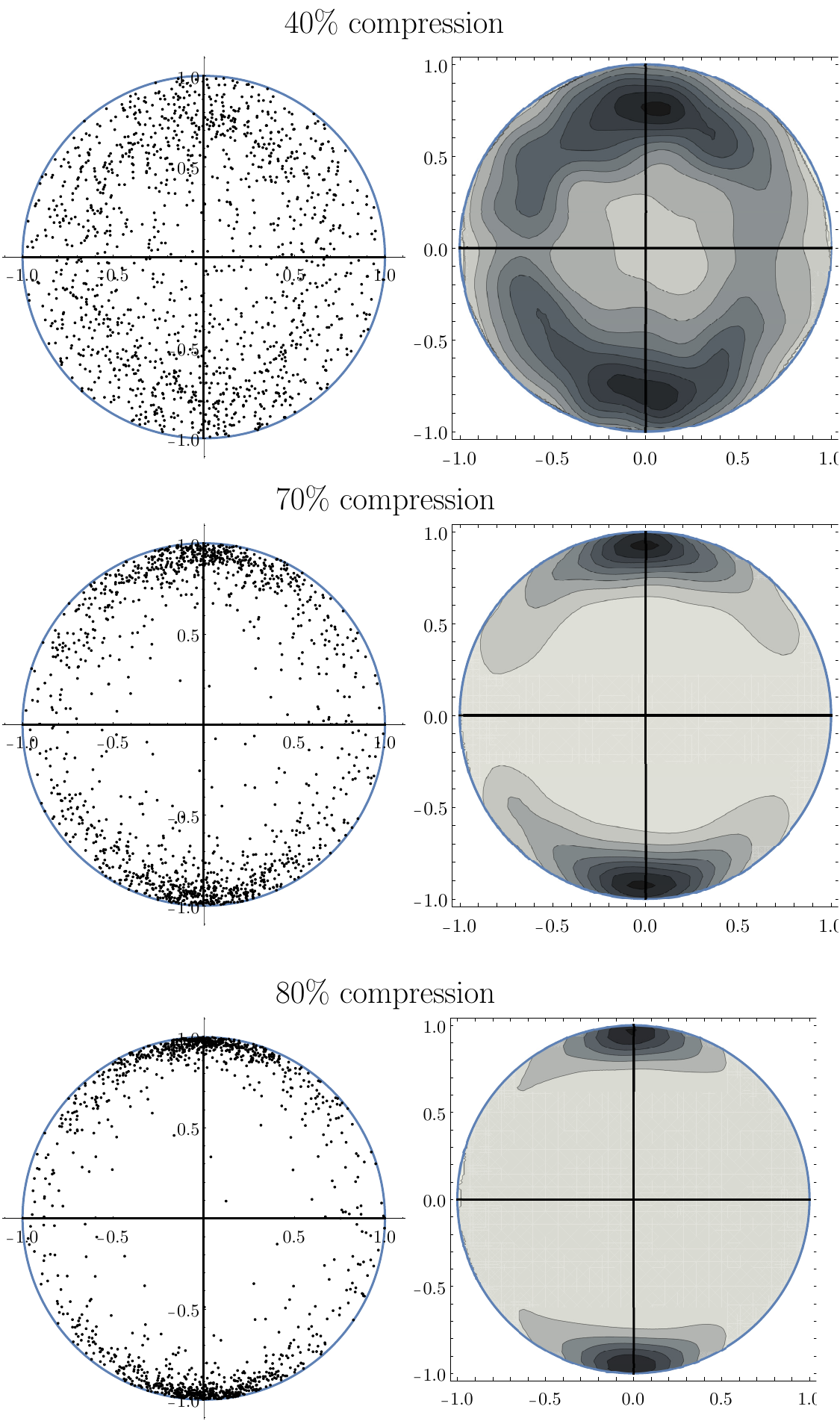}
\par\end{centering}
\caption{\label{figaTexture--planes}Texture $(111)$ planes, see also Alankar
\emph{et al.} \cite{alankar2009}. }

\end{figure}

\begin{figure}
\begin{centering}
\begin{tabular}{ccc}
\includegraphics[width=6cm]{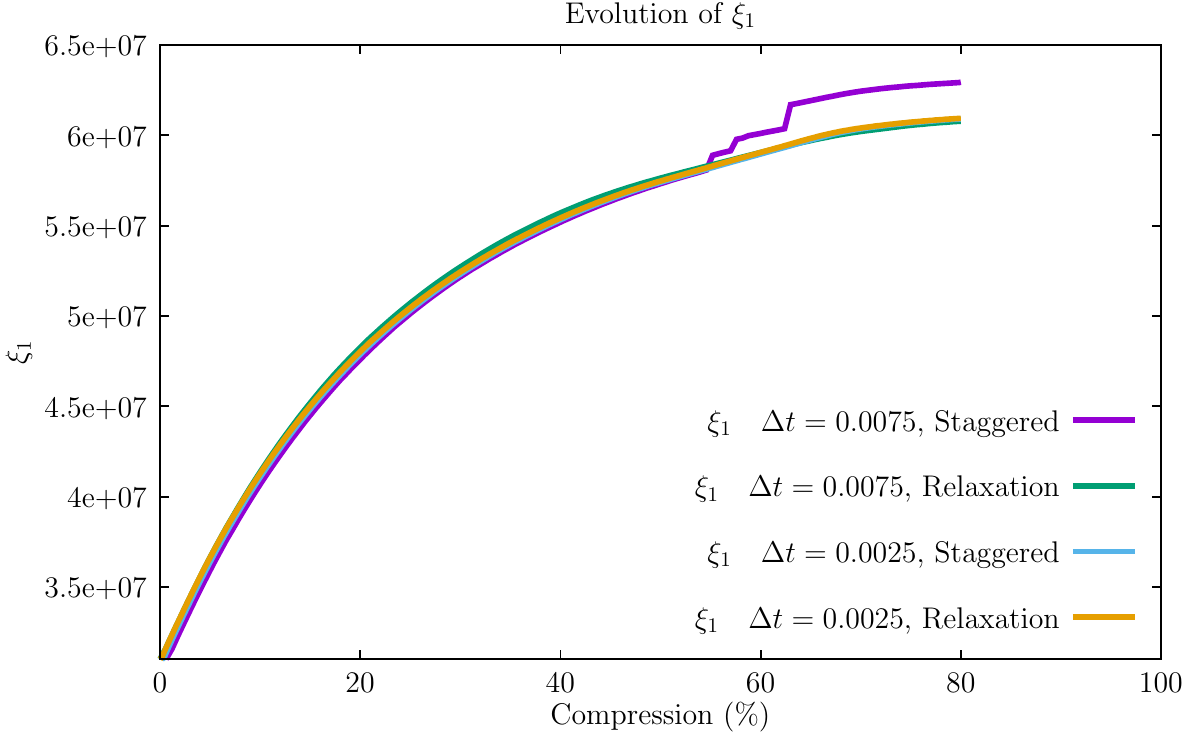} & \includegraphics[width=6cm]{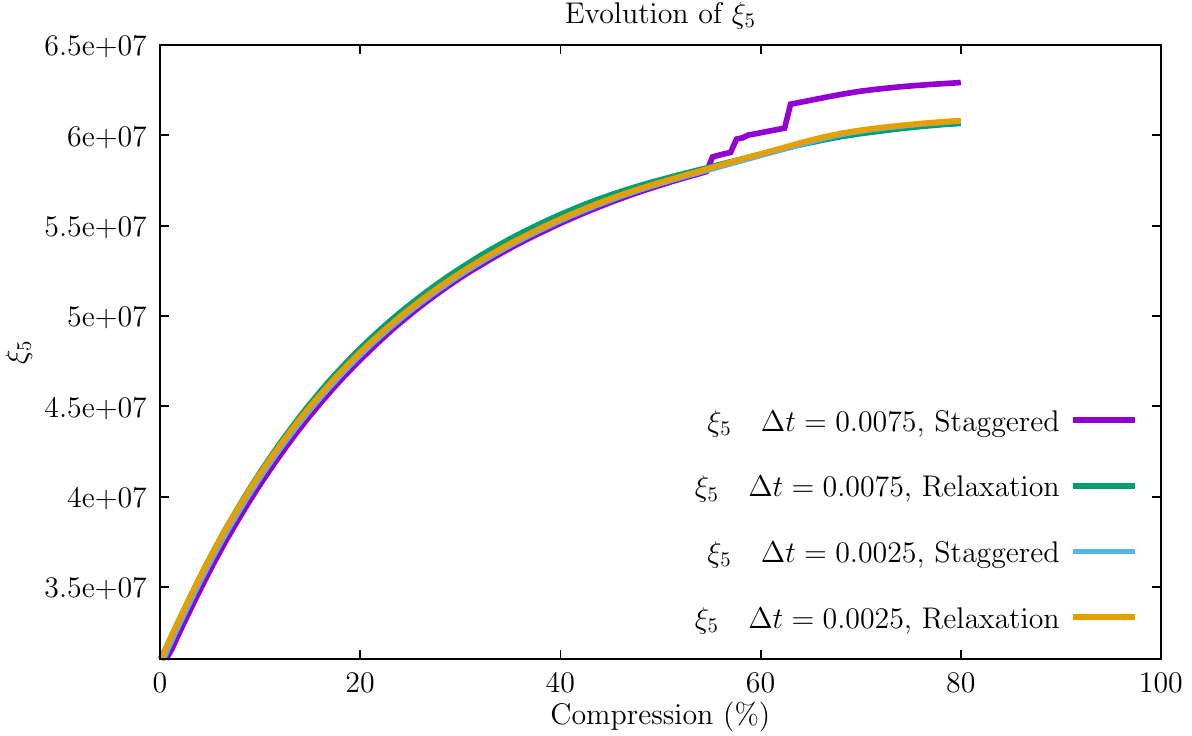} & \tabularnewline
\includegraphics[width=6cm]{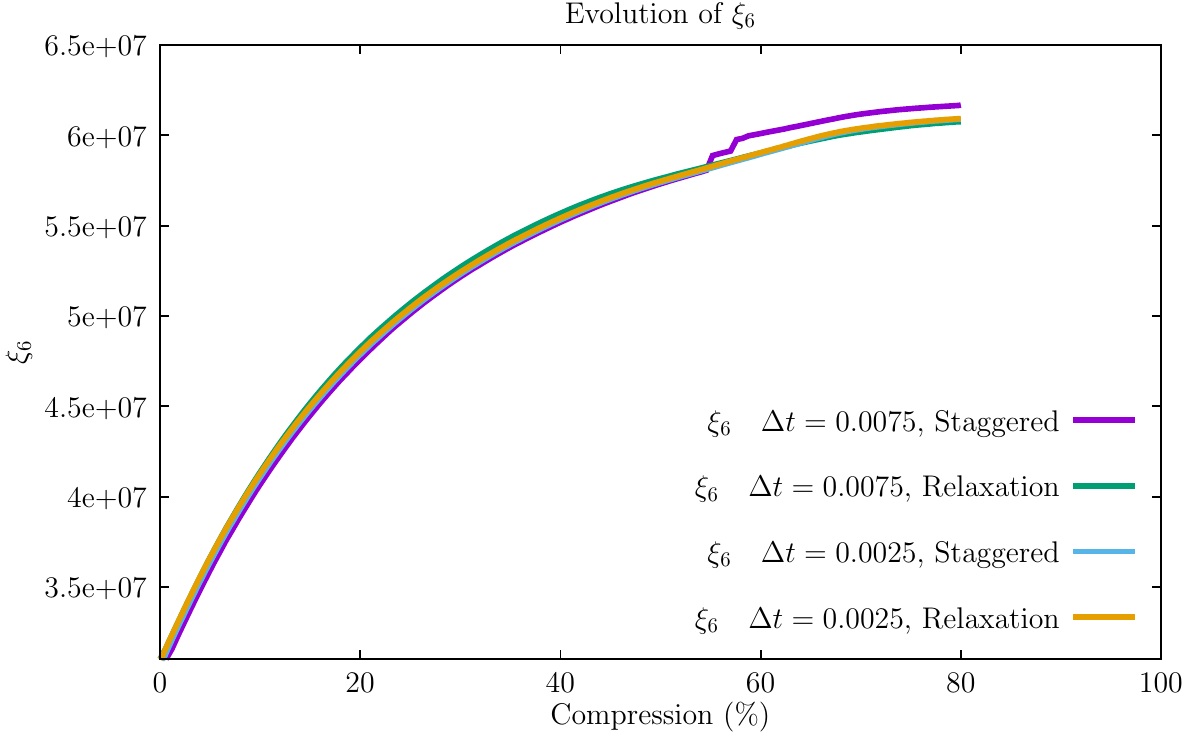} & \includegraphics[width=6cm]{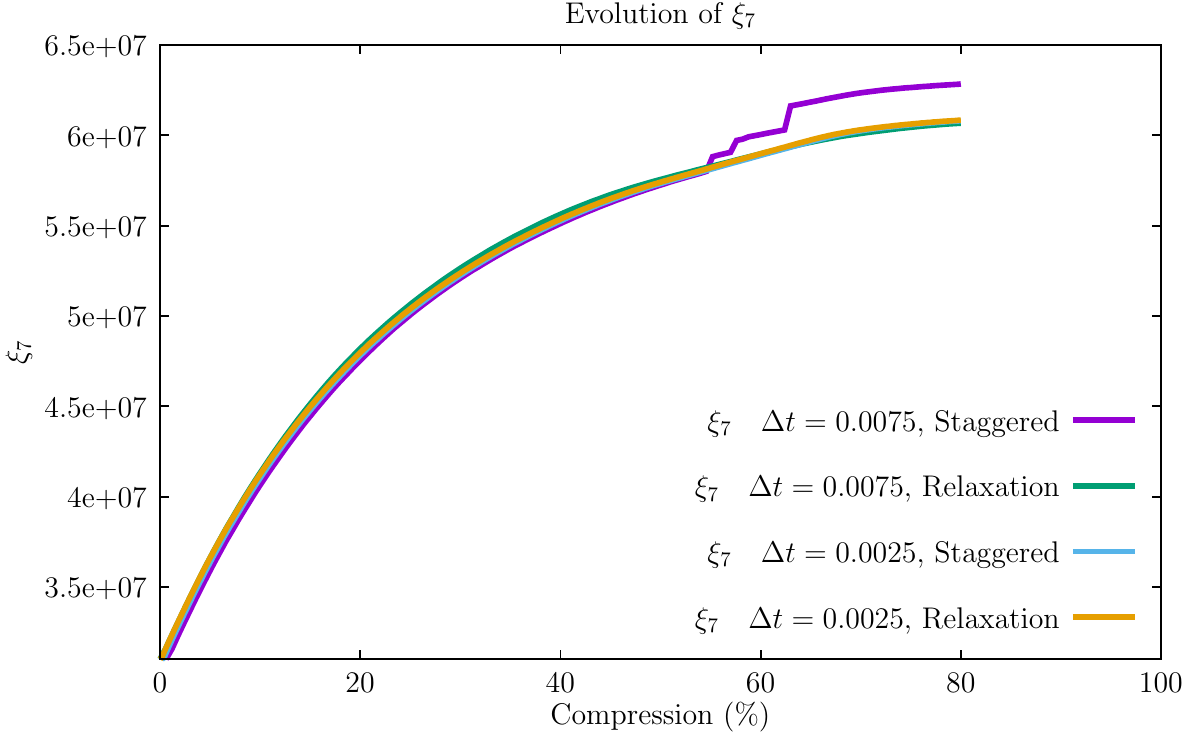} & \tabularnewline
\includegraphics[width=6cm]{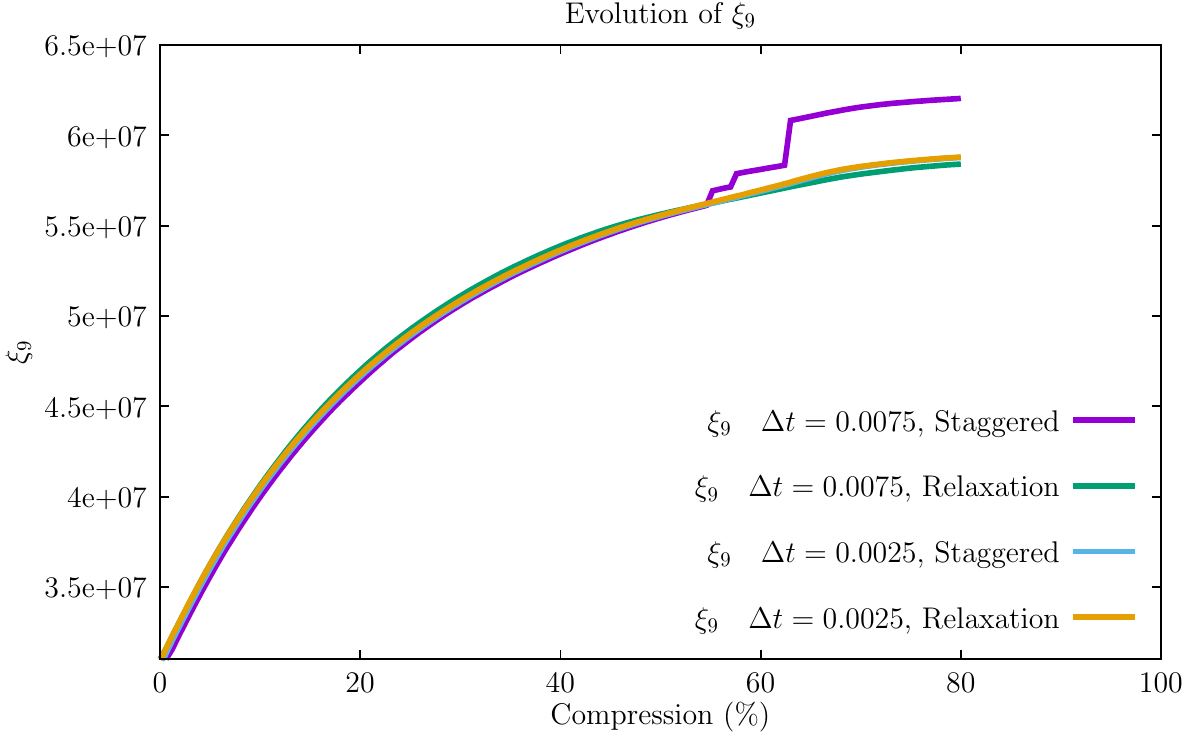} & \includegraphics[width=6cm]{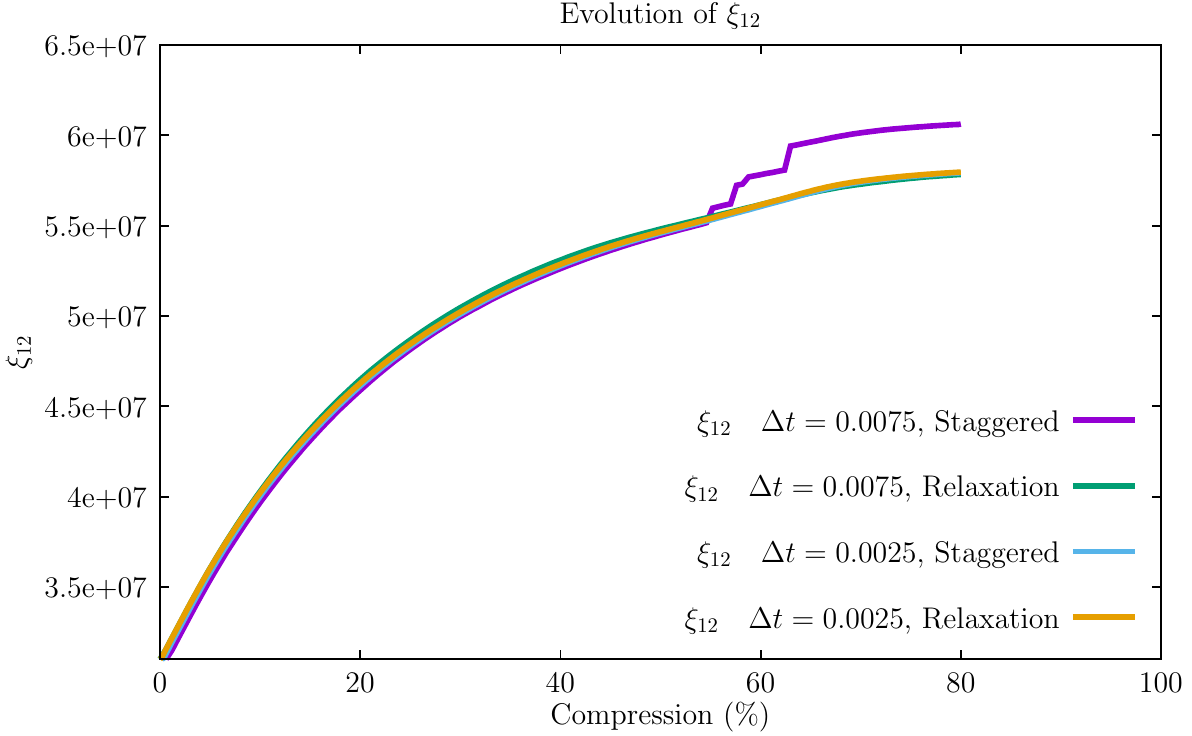} & \tabularnewline
\end{tabular}
\par\end{centering}
\caption{\label{figaEffect-of-the}Effect of the algorithm in the evolution
of $\xi_{1},\xi_{5,}\xi_{6},\xi_{7},\xi_{9}$ and $\xi_{12}$ in point
$X=2\times10^{-4}$, $Y=0$ and $Z=0$. }

\end{figure}

Texture evolution is shown in Figure \ref{figaTexture--planes} and
is similar to that published by Alankar \emph{et al.} \cite{alankar2009},
who used slightly different boundary conditions. We now assess drifting
of $\xi_{1},\cdots,\xi_{12}$. For conciseness reasons, Figure \ref{figaEffect-of-the}
shows the evolution of $\xi_{1},\xi_{5,}\xi_{6},\xi_{7},\xi_{9}$
and $\xi_{12}$ using the classical staggered algorithm and dynamic
relaxation. Two time-steps are tested: $\Delta t=0.0025$ s and $\Delta t=0.0075$
s. If the traditional staggered decomposition is adopted, with all
variables, for the larger time-step, at a compression value of $55\%$,
drifting occurs. Dynamic relaxation removes this effect and allows
the use of large time-steps. We also noted that Newton convergence
is improved with dynamic relaxation.

\clearpage

\section{Conclusions}

\label{seconclusion}

We introduced an alternative method to solve FCC single crystal finite
strain plasticity problems. It combines the use of Logarithmic strain
additive decomposition, Newton iteration for the plastic strain, a
linear solution for the evolution of the hardening variables and an
implicit staggered algorithm based on dynamic relaxation. Discretization
makes use of an EAS formulation based on the Green-Lagrange strain.
Verification tests for mesh size and time step dependences were successfully
performed and a polycrystal example from \cite{alankar2009} was studied
for texture evolution. We conclude that the dynamic relaxation is
effective in reducing drift caused by the staggered algorithm. Significant
savings and computational cost reductions can be achieved and the
procedure can be extended to more intricate constitutive laws.

\paragraph{Funding}
The authors acknowledge Funda\c{c}\~{a}o para a Ci\^{e}ncia e a Tecnologia (FCT)
for its financial support via the project LAETA Base Funding (DOI:
10.54499/UIDB/50022/2020)
\paragraph{Authors's contributions}
P. Areias carried out the bibliographic review and algorithm implementation, as well as the drafting. Charles dos Santos organized the bibliography by theme and relative contributions and drew the scientific diagrams. R. Melicio assisted with the writing and technical content and N. Silvestre performed the revisions, provided criticism and further references.
\paragraph{Competing interests}
\begin{itemize} 
\item The authors declare that the research was conducted in the absence of any commercial or financial relationships that could be construed as a potential conflict of interest.
\item The authors declare that they have no competing interests.
\end{itemize}
\paragraph{Open Access} This article is licensed under a Creative Commons Attribution 4.0 International License, which permits use, sharing, adaptation, distribution and reproduction in any medium or format, as long as you give appropriate credit to the original author(s) and the source, provide a link to the Creative Commons license, and indicate if changes were made. The images or other third party material in this article are included in the article's Creative Commons license, unless indicated otherwise in a credit line to the material. If material is not included in the article's Creative Commons license and your intended use is not permitted by statutory regulation or exceeds the permitted use, you will need to obtain permission directly from the copyright holder. To view a full copy of this license, visit http://creativecommons.org/licenses/by/4.0/.

 
\end{document}